\documentclass[A4paper]{article}

\usepackage{main}

\begin{document}

%\linenumbers

\title{Effective polygonal mesh generation and refinement for VEM}

\author{S. Berrone\footremember{trailer}{Dipartimento di Scienze Matematiche ``Giuseppe Luigi Lagrange'', Politecnico di Torino, Corso Duca degli Abruzzi 24, 10129 Torino, Italy (stefano.berrone@polito.it, fabio.vicini@polito.it).}, F. Vicini\footrecall{trailer}{}}

\maketitle

\begin{abstract}
%% Text of abstract
In the present work we introduce a novel refinement algorithm for two-dimensional elliptic partial differential equations discretized with Virtual Element Method (VEM).
The algorithm improves the numerical solution accuracy and the mesh quality through a controlled refinement strategy applied to the generic polygonal elements of the domain tessellation.
The numerical results show that the outlined strategy proves to be versatile and applicable to any two-dimensional problem where polygonal meshes offer advantages.
In particular, we focus on the simulation of flow in fractured media, specifically using the Discrete Fracture Network (DFN) model.
A residual a-posteriori error estimator tailored for the DFN case is employed.
We chose this particular application to emphasize the effectiveness of the algorithm in handling complex geometries. 
All the numerical tests demonstrate optimal convergence rates for all the tested VEM orders.

\end{abstract}


%% MSC codes here, in the form: \MSC code \sep code
%% or \MSC[2008] code \sep code (2000 is the default)
Mesh adaptivity, Virtual Element Method, Polygonal mesh refinement, Convergence and Optimality

\section{Introduction}
\label{sec:intro}
The necessity for the numerical resolution to partial differential equations (PDEs) in highly complex geometries arises from many applications in both science and engineering.
Dealing with these domains, characterized by intricate geometries and where the generation of a high-quality triangular mesh is either computationally expensive or impractical, has led to the adoption of general polygonal meshes in recent years.
In this context, there is a natural demand for a numerical method capable to tackle generic polygonal tessellation, and one of such method is the Virtual Element Method (VEM), \cite{doi:10.1142/S0218202516500160,doi:10.1142/S0218202512500492}.

Moreover, it is widely recognized that a posteriori error estimates are essential for ensuring the reliability of simulation tools.
Indeed, ad-hoc adaptive algorithms can be developed to reduce computational effort by leveraging the ability to control the approximation error with cost-effective computable quantities localized on mesh cells.

The study of adaptive schemes in conjunction with polygonal elements has garnered significant attention over the past few years.
In particular, in the VEM literature, numerous works have emerged following the introduction of the a-posteriori estimator in \cite{Cangiani2017}.
These works span both theoretical \cite{doi:10.1137/21M1458740, doi:10.3934/mine.2023101} and practical \cite{BERRONE2021103502} perspectives.

In this study, we present a novel refinement algorithm, designed for generic two-dimensional second order elliptic PDEs, extending the idea introduced in \cite{BERRONE2021103502} and \cite{BERRONE2022103770}.
Our approach starts with a generic polygonal mesh featuring the minimum number of Degrees of Freedoms (DOFs) consistent with the domain geometry, while meeting the minimum regularity requirements for the VEM convergence.
The refinement process involves the application of an optimal strategy to enhance both the mesh quality and the numerical solution accuracy.
This procedure is controlled by two parameters, allowing for a smooth transition of the tessellation from a poor quality complex polygonal mesh to a quasi-triangular or quasi-quadrilateral one within a few steps.

We consider the application of our proposed refinement algorithm to the simulation of the flow in fractured media.
The Discrete Fracture Network (DFN) model is employed.
We take advantage of the residual a-posteriori error estimator introduced in \cite{doi:10.1142/S0218202517500233} and extended to the DFN case in \cite{BERRONE2021103502}.
We remark that we choose this model problem due to its intrinsic generation of highly complex geometries and a large number of aligned edges in the initial tessellation, \cite{pichot2018simulations, doi:10.1137/18M1228736, doi:10.1007/s10092-023-00517-5}.
However, we stress that the strategy outlined in this paper is versatile and can be applied to any two-dimensional problem where polygonal meshes offer significant advantages. 

The outline of the paper unfolds as follows.
In Section~\ref{sec:setting}, we introduce the problem's context, the VEM discrete space, and the adaptive scheme.
In Section~\ref{sec:ref}, we focus on the newly proposed refinement strategy. 
Section~\ref{sec:problem} describes the DFN problem and the mesh generation on the media.
This section also introduces the a-posteriori error estimator.
Finally, in Section~\ref{sec:results}, we present the numerical results across three distinct cases of increasing complexity.
Namely a generic two-dimensional application and two scenarios involving DFNs.

Throughout this work, we use the following notations.
Given an open set $\omega \subset \mathbb{R}^2$ we indicate with $(\cdot, \cdot)_{\omega}$ and with $\left|\left|\cdot \right|\right|_{\omega}$ the $L^2(\omega)$ inner product and norm, respectively.

\section{Setting and Target}
\label{sec:setting}
We start considering, for this and the following Section~\ref{sec:ref}, a second order Poisson's problem defined on a generic planar polygonal domain $\Omega \subset \mathbb{R}^2$.    
We prescribe on the whole boundary, $\partial \Omega$, zero Dirichlet boundary conditions, hence we look for the variational solution $U$ living in the space $V := H^1_0$.
Different second order Partial Differential Equations (PDE) and boundary conditions can be considered. 

For the discrete approximation of $U$, we assume the existence of a tessellation $\mathcal{T}_{\Omega}$ on $\Omega$, composed by a finite number of polygonal cells $E \in \mathcal{T}_{\Omega}$.
We admit on $\mathcal{T}_{\Omega}$ the existence of aligned edges in the same element.

Throughout this section and the subsequent discussion, we use $\# \mathcal{T}_{\Omega}$ to indicate the total number of polygons in the mesh $\mathcal{T}_{\Omega}$.
We group all the edges of the tessellation into the set $\mathcal{E}_{\Omega}$.
The set of neighbouring cells of an edge $e \in \mathcal{E}_{\Omega}$ is denoted by $\mathcal{N}_e$.
It is worth noting that, in the application under consideration of Section~\ref{sec:problem}, the number of neighbouring cells $\# \mathcal{N}_e$ for and edge $e \in \mathcal{E}_{\Omega}$ might be higher than two (see Figure~\ref{fig:tipetut:network}).
Considering a generic polygon $E \in \mathcal{T}_{\Omega}$, we use $N_E$ to denote its number of vertices and edges.
The set of all the edges of $E$ is indicated by $\mathcal{E}_E$.
Additionally, $H_E$ and $h_E$ refer to the longest and the smallest edge length $|e|$, respectively, for all $e \in \mathcal{E}_E$.
When selecting an edge $e \in \mathcal{E}_E$, $\mathcal{I}_E^e$ represents the collection of all the edges of $E$ aligned with $e$.
Moreover, $|\mathcal{I}_E^e|$ stands for the sum of the lengths of these aligned edges.
The centroid of the element $E$ is denoted by $\mathbf{x}_E$, and $\mathbb{I}_E$ indicates the inertia tensor with respect to $\mathbf{x}_E$.
We account by $r_E$ for the minimum distance between the centroid $\mathbf{x}_E$ and the edges of $E$ and by $R_E$ the maximum distance between the centroid $\mathbf{x}_E$ and the vertices of $E$.
Finally, $D_E$ refers to the diameter of the element $E$ and we define the \emph{mesh size} of $\mathcal{T}_{\Omega}$ as $D := \max_{E \in \mathcal{T}_{\Omega}} D_E$.

\subsection{VEM Space}
We introduce the Virtual Element Method (VEM) to address the computation of a discrete solution on a polygonal mesh $\mathcal{T}_{\Omega}$.
For a given $k \in \mathbb{N}^+$, we denote by $\mathbb{P}_k(E)$ the set of polynomials of degree up to $k$ defined on each element $E \in \mathcal{T}_{\Omega}$.
With the same notation used in \cite{AHMAD2013376}, we define the projectors $\Pi^{\nabla}_k : H^1(E) \to \mathbb{P}_k(E)$ and $\Pi^0_k : L^2(E) \to \mathbb{P}_k(E)$ such that
\begin{linenomath}
    \begin{equation*}
        \begin{cases}
            \left(\nabla \Pi^{\nabla}_k v - \nabla v, \nabla p  \right)_E = 0 & \forall p \in \mathbb{P}_k,\\
            (\Pi^{\nabla}_k v - v, 1)_{\partial E} = 0 & \text{if } k = 1,\\
            (\Pi^{\nabla}_k v - v, 1)_E = 0 & \text{if } k > 1,
        \end{cases}
    \end{equation*}
\end{linenomath}
and
\begin{linenomath}
    \begin{equation*}
        \left(\Pi^0_k v - v, p  \right)_E = 0 \ \forall p \in \mathbb{P}_k.
    \end{equation*}
\end{linenomath}
As done in \cite{doi:10.1142/S0218202516500160}, $\forall E \in \mathcal{T}_{\Omega}$ we introduce the local VEM space as 
\begin{linenomath}
    \begin{equation}
        \label{eq:space:discrete:local}
        V^E_k = \left\{v \in H^1(E) : v_{|_e} \in \mathbb{P}_k(e) \ \forall e \in \mathcal{E}_E,\ v_{|_{\partial E}} \in C^0(\partial E),\ \Delta v \in \mathbb{P}_k(E),\ \left(v-  \Pi^{\nabla}_k v, p \right)_E = 0 \ \forall p \in \mathbb{P}_k(E) / \mathbb{P}_{k-2}(E)  \right\},
    \end{equation}
\end{linenomath}
where $\mathbb{P}_k(E) / \mathbb{P}_{k-2}(E)$ denotes the polynomials $\mathbb{P}_k(E)$ that are $L^2$-orthogonal to $\mathbb{P}_{k-2}(E)$.

Finally, we define the global discrete space $V_k \subseteq V$ as:
\begin{linenomath}
    \begin{equation}
        \label{eq:space:discrete}
        V_k = \left\{ v \in V : v \in C^0(\Omega),\ v_{|_E} \in V^E_k,\ \forall E \in \mathcal{T}_{\Omega} \right\}.
    \end{equation}
\end{linenomath}
We uniquely identify a function $v \in V_k$ using the set of Degrees of Freedoms (DOFs):
\begin{itemize}
    \item the values of $v$ at $N_E$ vertices of $E$;
    \item if $k > 1$, the values of $v$ at $k - 1$ Gauss internal quadrature points of each of the $N_E$ edges of $E$;
    \item if $k > 1$, the internal scaled moments of order $k-2$ of $v$ in $E$.
\end{itemize}
This choice is not unique, see an other example in \cite{doi:10.1142/S0218202516500160}, however the selected DOFs are unisolvent for $V_k$ \cite{doi:10.1142/S0218202512500492}.

\subsection{Target}
The goal of this work is to present an adaptive procedure based on the well-known scheme
\begin{linenomath}
    \begin{equation}
        \label{eq:schema}
        \texttt{SOLVE} \to \texttt{ESTIMATE} \to \texttt{MARK} \to \texttt{REFINE}.
    \end{equation}
\end{linenomath}
The procedure we propose is tailored for a generic problem discretized using the VEM schema, but can be extended to different polygonal methods. 
We start the Process~\eqref{eq:schema} with a generic polygonal mesh $\mathcal{T}^0_{\Omega}$ which presents the number of DOFs as low as possible.
To achieve this target, we opt for $\mathcal{T}^0_{\Omega}$ as the convex polygonal approximation of $\Omega$.
This choice minimizes the initial number of DOFs, as suggested in \cite{BERRONE2022103770}.
In what follows, we call this particular choice \emph{minimal mesh}.
It is worth noting that if $\Omega$ is convex, then $\mathcal{T}^0_{\Omega}$ can be $\Omega$ itself.

The modules of Process~\eqref{eq:schema} are defined as follows: for each iteration $m \geq 0$, given a mesh $\mathcal{T}^m_{\Omega}$ 
\begin{itemize}
    \item $[u^m] = \texttt{SOLVE}(\mathcal{T}_{\Omega}^m)$ creates the VEM discrete solution $u^m \in V_k$ of the PDE;
    \item $[\eta_{\Omega}^m] = \texttt{ESTIMATE}(\mathcal{T}_{\Omega}^m, u^m)$ computes an estimation $\eta^m$ of the error between $U$ and $u^m$ ;
    \item $[\mathcal{M}^m] = \texttt{MARK}(\mathcal{T}_{\Omega}^m, \eta_{\Omega}^m)$ selects a subset of cells $\mathcal{M}^m \subseteq \mathcal{T}_{\Omega}^m$ candidates for refining;
    \item $[\mathcal{T}_{\Omega}^{m + 1}] = \texttt{REFINE}(\mathcal{T}_{\Omega}^m, \mathcal{M}^m)$ creates a new fine-tuned discretization splitting the marked elements.
\end{itemize}

We aim the \texttt{REFINE} rule to be applied to generic (convex) polygonal elements and to preserve or improve the mesh quality.
In cases the initial mesh $\mathcal{T}^0_{\Omega}$ exhibits poor quality, we aim to enhance this quality during subsequent refinement process iterations.
It is crucial to note that mesh quality is essential for producing a robust VEM solution in the \texttt{SOLVE} procedure.

We require that the initial mesh $\mathcal{T}^0_{\Omega}$ consists of convex elements and satisfies the minimal quality conditions essential for the VEM convergence, as outlined in \cite{doi:10.1007/s10444-021-09913-3}.
Throughout this discussion, we omit the superscript $m$ denoting the iteration in the mesh.

The literature indicates that the VEM convergence is guaranteed combining the condition:
\begin{enumerate}
    \item there exists $\gamma_r \in (0, 1)$ independent of $D$ such that
    \begin{equation}
        \tag{C1}\label{eq:minimal:conditions:r}
        r_E \geq \gamma_r D_E, \forall E \in \mathcal{T}_{\Omega},
    \end{equation}
\end{enumerate}
with either one of the conditions:
\begin{enumerate}
    \setcounter{enumi}{1}
    \item there exists $\gamma_h \in (0, 1)$ independent of $D$  such that
    \begin{equation}
        \tag{C2}\label{eq:minimal:conditions:h}
        h_E \geq \gamma_h D_E, \forall E \in \mathcal{T}_{\Omega};
    \end{equation}
    \item there exists $N \in \mathbb{N}^+$ independent of $D$ such that
    \begin{equation}
        \tag{C3}\label{eq:minimal:conditions:n}
        N_E \leq N, \forall E \in \mathcal{T}_{\Omega};
    \end{equation}
    \item there exists $\gamma_{al} \in \mathbb{R}^+$ independent of $D$ such that $\forall E \in \mathcal{T}_{\Omega}$ the boundary $\partial E$ can be sub-divided in $K$ finite disjoint sequence of edges $\mathcal{I}_E^1, \dots, \mathcal{I}_E^K$ such that
    \begin{equation}
        \tag{C4}\label{eq:minimal:conditions:al}
        \frac{\max_{e \in \mathcal{I}_E^k} |e|}{\min_{e \in \mathcal{I}_E^k} |e|} \leq \gamma_{al}, \forall k \in \{1,\dots, K\}.
    \end{equation}
\end{enumerate}
Condition~\eqref{eq:minimal:conditions:al} can be interpreted in term of the \emph{piecewise quasi-uniformity} of consecutive edges $\mathcal{I}_E^e$ aligned to $e \in \mathcal{E}_E$ (refer to \cite{doi:10.1137/21M1411275} for more details). 

We stress that the convexity of the tessellation elements and conditions~\eqref{eq:minimal:conditions:r}-\eqref{eq:minimal:conditions:al} do not inherently ensure high-quality mesh elements.
Indeed, the elements of $\mathcal{T}_{\Omega}$ may exhibit undesirable \emph{badly-shapes} characteristics, such as collapsing bulks, small edges, an excess of edges, or irregular aligned edge lengths.
Badly shaped elements are characterized by the constants $\gamma_r$ and $\gamma_h$ being close to zero or to the values of the constants $N$, and $\gamma_{al}$ large.

To address this, we propose a new \texttt{REFINE} strategy with the following objectives:
\begin{itemize}
    \item Generate new sub-polygons with improved quality in term of $\gamma_r$ and $\gamma_h$, if low;
    \item Control the number of vertices in the sub-polygons, aiming for their reduction during the refinement iterations;
    \item Limit and decrease the number of aligned edges, or at least produce quasi-uniform sets $\mathcal{I}^k_E$, $\forall E \in \mathcal{T}^m_{\Omega}$.
\end{itemize}
\section{The \texttt{REFINE} module}
\label{sec:ref}
\begin{algorithm}[!h]
    \caption{$[\mathcal{T}_{\Omega}^{m+1}] = \texttt{REFINE}(\mathcal{T}_{\Omega}^m, \mathcal{M}^m, c_{\rho}, c_{al})$}\label{alg:ref}
    \begin{algorithmic}[1]
        \REQUIRE the mesh $\mathcal{T}_{\Omega}^m$, the marked elements $\mathcal{M}^m$, the parameters $\{c_{\rho}, c_{al}\}$
        \ENSURE the refined mesh $\mathcal{T}_{\Omega}^{m+1}$
        \FOR {$E \in \mathcal{M}^m$}
        \STATE Compute $\hat{E}$ as the polygon formed by the union of aligned edges of $E$
        \STATE Compute the \texttt{MAX-MOMENTUM}$(\hat{E})$ direction $\mathbf{t}_{mm}$
        \STATE Fix the direction $\mathbf{t}$ using \texttt{SMOOTH-DIRECTION}$(E, \mathbf{t}_{mm}, c_{\rho})$
        \STATE Split $E = E_1 \cup E_2$ with the direction $\mathbf{t}$, deactivate $E$
        \STATE Mark the edges $e \in \mathcal{E}_{E_1} \cup \mathcal{E}_{E_2}$ split by $\mathbf{t}$
        \STATE Update the polygon neighbours $\{\mathcal{N}_e \setminus \{E_1, E_2\}, \forall e \in \mathcal{E}_{E_1} \cup \mathcal{E}_{E_2} \text{ marked }\}$
        \STATE Set $\mathcal{Q} = \{\mathcal{N}_e, \forall e \in \mathcal{E}_{E_1} \cup \mathcal{E}_{E_2} \text{ marked }\}$
        \FOR {$E_q \in \mathcal{Q}$}
        \FOR {$e \in \mathcal{E}_{E_q}$ marked}
        \IF {\texttt{CHECK-QUALITY}$(E_q, e, 1, c_{\rho}, c_{al})$ failed} 
        \STATE Set $\mathcal{M}^m = \mathcal{M}^m \cup E_q$ \label{alg:ref:extend}
        \ELSE
        \STATE Unmark $e$
        \ENDIF
        \ENDFOR
        \ENDFOR
        \ENDFOR
    \end{algorithmic}
\end{algorithm}
\begin{algorithm}[!h]
    \caption{$[\mathbf{t}_{mm}] = \texttt{MAX-MOMENTUM}(\hat{E})$}\label{alg:max_mom}
    \begin{algorithmic}[1]
        \REQUIRE the polygon $\hat{E}$ with no aligned edges
        \ENSURE the cut direction $\mathbf{t}_{mm}$
        \IF {$\hat{E}$ is triangle}
        \STATE Return $\mathbf{t}_{mm}$ as the \texttt{NVB}$(\hat{E})$ direction
        \ENDIF
        \STATE Compute the inertia tensor $\mathbb{I}_{\hat{E}}$ respect the centroid $\mathbf{x}_{\hat{E}}$
        \STATE Set $\mathbf{t}_{mm}$ the line passing by $\mathbf{x}_{\hat{E}}$ and parallel to the eigenvector of the maximum eigenvalue of $\mathbb{I}_{\hat{E}}$
    \end{algorithmic}
\end{algorithm}
\begin{algorithm}[!h]
    \caption{$[\mathbf{t}] = \texttt{SMOOTH-DIRECTION}(E, \mathbf{t}_{cut}, c_{\rho}, c_{al})$}\label{alg:smooth}
    \begin{algorithmic}[1]
        \REQUIRE the mesh polygon $E$, the cut direction $\mathbf{t}_{cut}$, the parameters $\{c_{\rho}, c_{al}\}$
        \ENSURE the smoothed direction $\mathbf{t}$
        \IF {$E$ is triangle}
        \STATE Return $\mathbf{t}_{cut}$ 
        \ENDIF
        \STATE Set $\mathbf{t} = \mathbf{t}_{cut}$
        \FOR {$e \in \mathcal{E}_E$ intersected by $\mathbf{t}_{cut}$}
        \IF {\texttt{CHECK-QUALITY}$(E, e, 2, c_{\rho}, c_{al})$ failed}
        \STATE Move $\mathbf{t}$ to intersect the closest vertex of $e$
        \ELSE
        \STATE Move $\mathbf{t}$ to intersect the middle point of $e$
        \ENDIF
        \ENDFOR
    \end{algorithmic}
\end{algorithm}
\begin{algorithm}[!h]
    \caption{$[bool] = \texttt{CHECK-QUALITY}(E, e, s, c_{\rho}, c_{al})$}\label{alg:quality}
    \begin{algorithmic}[1]
        \REQUIRE the mesh polygon $E$, the edge $e \in \mathcal{E}_E$ to check, $s \in \mathbb{N}^+$ the number of $e$ uniform subdivisions, the parameters $\{c_{\rho}, c_{al}\}$
        \ENSURE True if $e$ quality is respected
        \STATE Set $\rho_e = \max_{E_n} \min \{h_{E_n}, r_{E_n}\}$, $\forall E_n \in \mathcal{N}_e$
        \IF {$|e| < c_{\rho} \rho_e s$ }
        \STATE Return False
        \ENDIF
        \STATE Set $\mathcal{I}_{E_n}^e$ the set of edges  $\tilde{e} \in \mathcal{E}_{E_n}$ aligned to $e$, $\forall E_n \in \mathcal{N}_e $
        \STATE Compute $|\mathcal{I}_{E_n}^e| = \sum_{e_a \in \mathcal{I}_{E_n}^e} |e_a|$ and $n_{E_n} = \# \mathcal{I}_{E_n}^e$
        \STATE Set $\mathcal{I}^s_e = \max_{E_n} |\mathcal{I}_{E_n}^e| / (n_{E_n} + (s - 1))$
        \IF {$|e| < c_{al}\mathcal{I}^s_e s$} \label{alg:ref:al}
        \STATE Return False
        \ENDIF
        \STATE Return True
    \end{algorithmic}
\end{algorithm}
We now discuss the primary contribution of this paper, i.e., the \texttt{REFINE} procedure.
Algorithm~\ref{alg:ref} outlines its pseudo-code.
We are going to provide a comprehensive description of the function and subsequently delve into the details of its sub-routines.

At the step $m$ of Process~\eqref{eq:schema}, the \texttt{REFINE} method starts from the mesh $\mathcal{T}_{\Omega}^m$, and the set of marked cells $\mathcal{M}^m$ obtained from the \texttt{MARK} module. 
For every marked cells $E \in \mathcal{M}^m$, the algorithm splits $E$ to decrease the local a-posteriori error estimator.
The objective is to bisect $E = \{E_1, E_2\}$ along a direction $\mathbf{t}$, determined by the combination of the outcomes of the \texttt{MAX-MOMENT} Algorithm~\ref{alg:max_mom} and \texttt{SMOOTH-DIRECTION} Algorithm~\ref{alg:smooth}.
Given the direction $\mathbf{t}$, we split the intersected edge of $\mathcal{E}_E$ producing new edges $e$ that are marked, and we update all the neighbouring cells $\mathcal{N}_e \setminus \{E_1, E_2\}$ of these marked edges $e$.

The VEM's capability to handle aligned edges simplifies the update of each cell in $\mathcal{N}_e \setminus \{E_1, E_2\}$, resulting in the creation of a new polygon containing the new aligned edges $e$ in place of the original one split by $\mathbf{t}$.
Subsequently, a set $\mathcal{Q}$ is generated, encompassing the neighbour cells $\mathcal{N}_e$ of the marked edges $e$.
For all the marked edges in the boundary of cells $E_q \in \mathcal{Q}$, the \texttt{CHECK-QUALITY} Algorithm~\ref{alg:quality} assesses the quality of the marked edges.
If the outcome is positive the edge is unmarked.
On the other hand, if the outcome is negative, the cell $E_q \in \mathcal{Q}$ is added to the set of the marked cells $\mathcal{M}^m$.
This operation, referred to hereafter as \emph{extension}, along with its incorporation into the Process~\eqref{eq:schema}, constitutes a key innovation in this work to ensure the best rate of convergence of the method.

For a generic polygon $\hat{E}$ devoid of aligned edges, the \texttt{MAX-MOMENT} function calculates the \emph{max-momentum} direction $\mathbf{t}_{mm}$.
This direction passes through the centroid $\mathbf{x}_{\hat{E}}$ of the element and is parallel to the eigenvector associated with the maximum eigenvalue of the inertia tensor $\mathbb{I}_E$.
An exception is made if the selected cell $\hat{E}$ is a triangle.
In this case, we employ the newest-vertex bisection (NVB) split criterion, see \cite{Nochetto2012}.
As demonstrated in \cite{BERRONE2022103770, doi:10.1142/S0218202517500233}, the \emph{max-momentum} direction has often the ability to generate sub-cells with improved quality compared to the original $\hat{E}$ for elongated elements, as measured by Condition~\eqref{eq:minimal:conditions:r}. 

The \texttt{SMOOTH-DIRECTION} Algorithm~\ref{alg:smooth} introduces a slight modification to the original \emph{max-momentum} direction $\mathbf{t}_{mm}$ computed on $\hat{E}$.
This adjustment involves shifting the cut direction to the middle of the edge $e \in \mathcal{E}_E$ intersected by $\mathbf{t}_{mm}$, or collapsing it to the closest vertex of $e$ if the \texttt{CHECK-QUALITY} Algorithm~\ref{alg:quality} fails.
As proved in \cite{BERRONE2022103770}, this strategy results in the generation of new polygon cells with, at most, the number of vertices of the original cell $E$.
The number of vertices only rises when the modified cut direction $\mathbf{t}$ splits two consecutive edges; in all the other cases, the children cells have fewer or an equal number of vertices as $E$, see \cite{BERRONE2022103770}.
This ensures control and improvement of Condition~\eqref{eq:minimal:conditions:n}.

The \texttt{CHECK-QUALITY} Algorithm~\ref{alg:quality} is applied to the polygonal cell $E$ targeted for splitting.
It assesses whether the edge $e \in \mathcal{E}_E$ can be divided in $s$ uniform parts.
In the \texttt{REFINE} algorithm we select $s=2$. Thus, the selected edge $e \in \mathcal{E}_E$ is split if the two conditions are satisfied:
\begin{eqnarray}
    \frac{|e|}{2} &\geq& c_{\rho} \max_{E_n \in \mathcal{N}_e} \min{\{h_{E_n}, r_{E_n}\}}\label{check:one},\\
    \frac{|e|}{2} &\geq& c_{al} \max_{E_n \in \mathcal{N}_e} \frac{|\mathcal{I}_{E_n}^e|}{\# \mathcal{I}_{E_n}^e + 1}\label{check:two}.
\end{eqnarray}
We recall that, $\mathcal{I}_E^e$ denotes the set of edges of $E$ aligned to $e$, and $|\mathcal{I}_E^e|$ and $\# \mathcal{I}_E^e$ represent the total length and the number of the contiguous aligned edge of $e$, respectively.

Check~\eqref{check:one} ensures that the splitting of the edge contributes to an improvement of the Condition~\eqref{eq:minimal:conditions:r} and Condition~\eqref{eq:minimal:conditions:h}.
The real parameter $c_{\rho} \geq 0.0$ is introduced to relax ($0.0 \leq c_{\rho} < 1.0$) or tighten  ($c_{\rho} \geq 1.0$) the constraint.
In addition, we propose Check~\eqref{check:two} to satisfy Condition~\eqref{eq:minimal:conditions:al}.
This check examines whether the newly created edges are \emph{quasi-uniform} compared to the other aligned edges.
The real parameter $c_{al} \geq 0.0$ controls the number of aligned edges permitted in the new refined mesh.
Namely, a value of $c_{al}$ lower than $1.0$ allows aligned edges, while $c_{al} > 1.0$ often rejects any aligned edges.

\subsection{Contributions of the new algorithm}
\begin{figure}[!h]
    \centering
    \begin{subfigure}[!h]{0.32\textwidth}
        \centering
        \includegraphics[width=\textwidth]{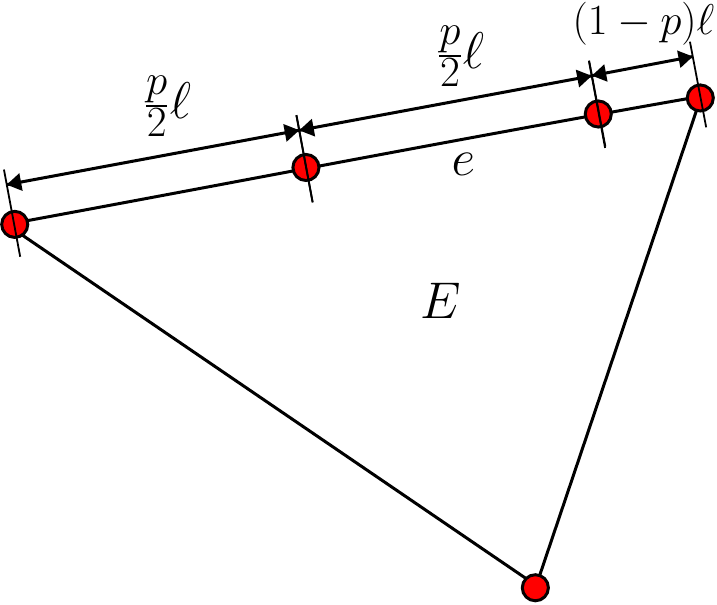}
        \caption{Original $E$}
        \label{fig:algo:aligned:original}
    \end{subfigure}
    \begin{subfigure}[!h]{0.32\textwidth}
        \centering
        \includegraphics[width=\textwidth]{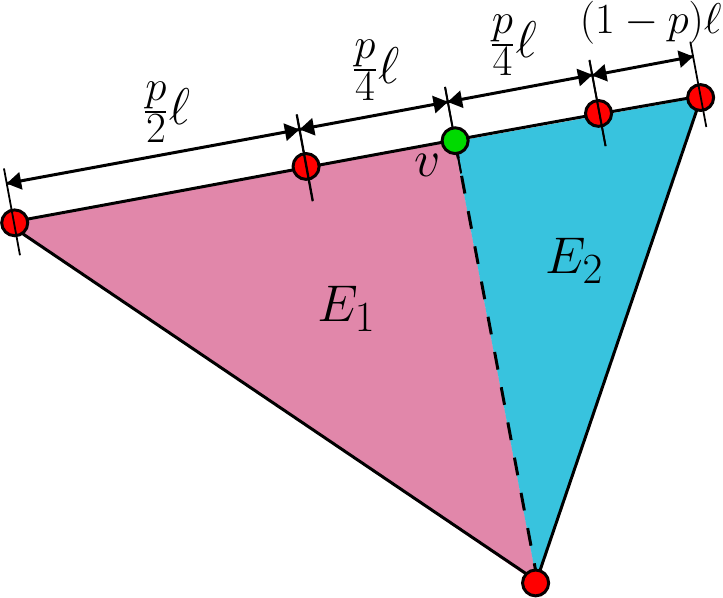}
        \caption{\emph{Old} algorithm}
        \label{fig:algo:aligned:old}
    \end{subfigure}
    \begin{subfigure}[!h]{0.32\textwidth}
        \centering
        \includegraphics[width=\textwidth]{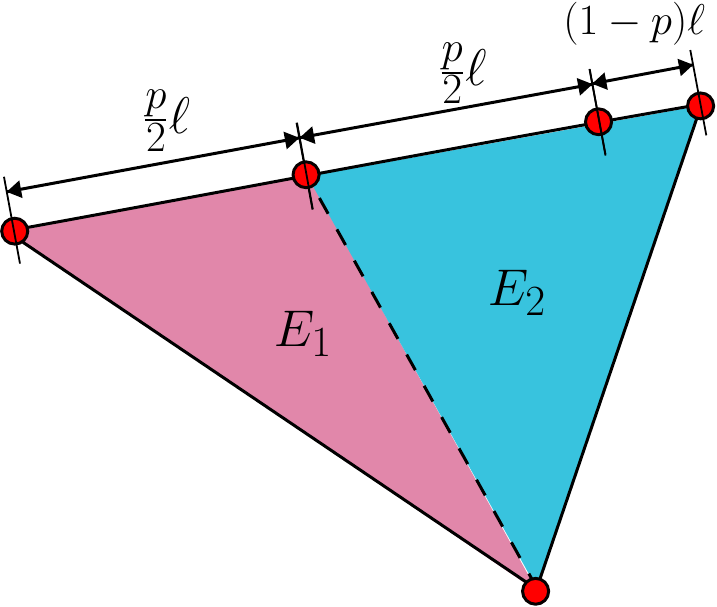}
        \caption{\emph{New} algorithm}
        \label{fig:algo:aligned:new}
    \end{subfigure}
    \caption{Refinement of $E$ - Case $1$ - $c_{\rho} = 0.5$ and $\frac{4}{5} < p < 1$}
    \label{fig:algo:aligned}
\end{figure}
\begin{figure}[!h]
    \centering
    \begin{subfigure}[!h]{0.32\textwidth}
        \centering
        \includegraphics[width=\textwidth]{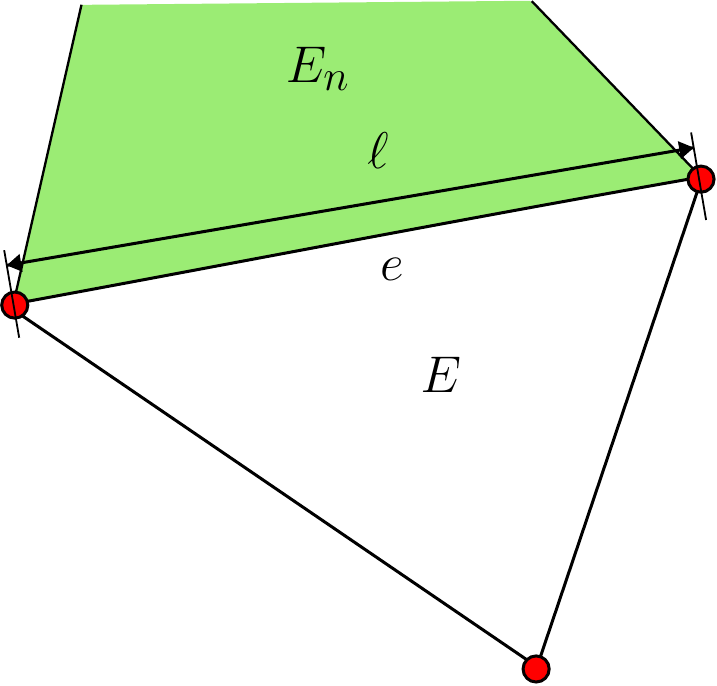}
        \caption{Original $E$}
        \label{fig:algo:neigh:original}
    \end{subfigure}
    \begin{subfigure}[!h]{0.32\textwidth}
        \centering
        \includegraphics[width=\textwidth]{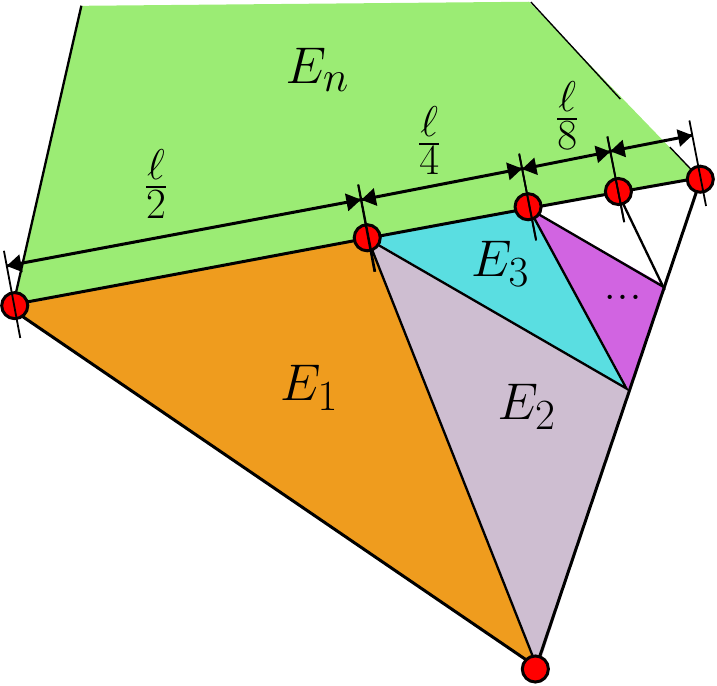}
        \caption{\emph{Old} algorithm}
        \label{fig:algo:neigh:old}
    \end{subfigure}
    \begin{subfigure}[!h]{0.32\textwidth}
        \centering
        \includegraphics[width=\textwidth]{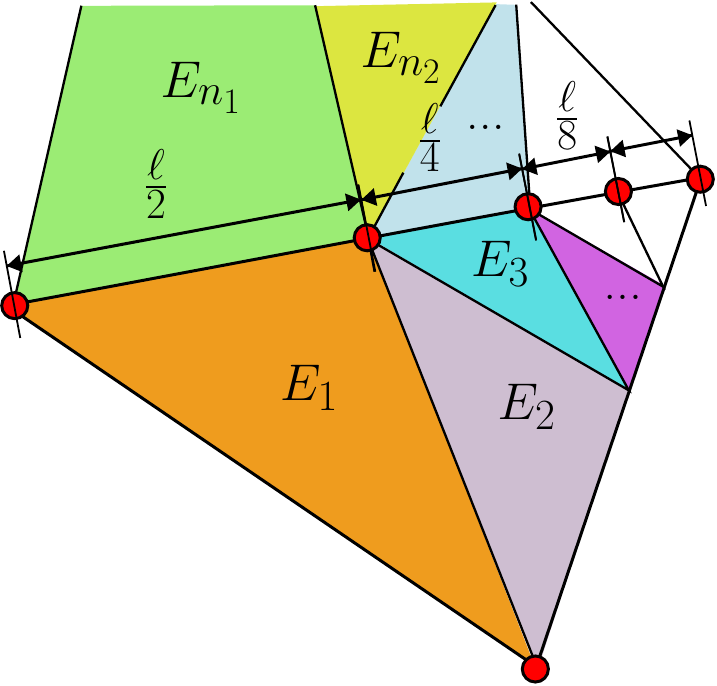}
        \caption{\emph{New} algorithm}
        \label{fig:algo:neigh:new}
    \end{subfigure}
    \caption{Refinement of $E$ - Case $2$ - $c_{\rho} = 0.5$}
    \label{fig:algo:neigh}
\end{figure}
\begin{figure}[!h]
    \centering
    \begin{subfigure}[!h]{0.32\textwidth}
        \centering
        \includegraphics[width=\textwidth]{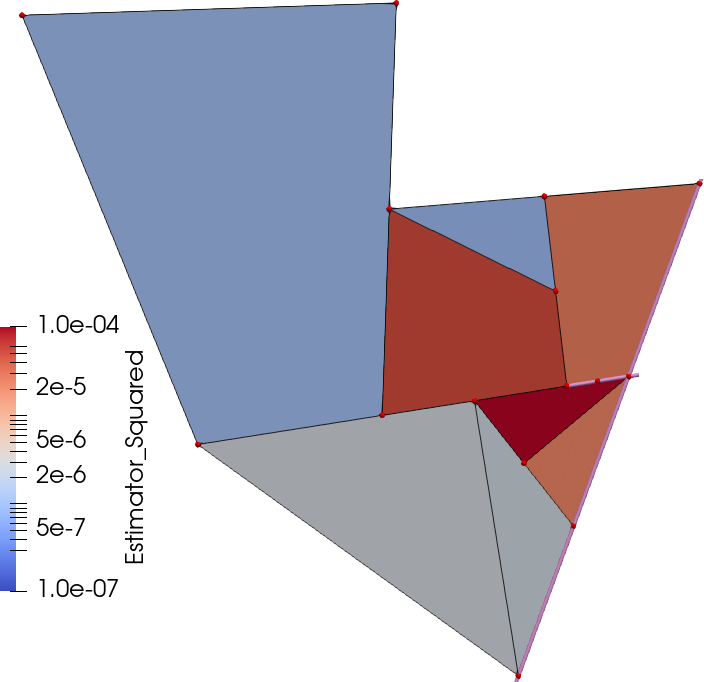}
        \caption{Case $1$ - \emph{Old} algorithm}
        \label{fig:algo:aligned:example:old}
    \end{subfigure}
    \begin{subfigure}[!h]{0.32\textwidth}
        \centering
        \includegraphics[width=\textwidth]{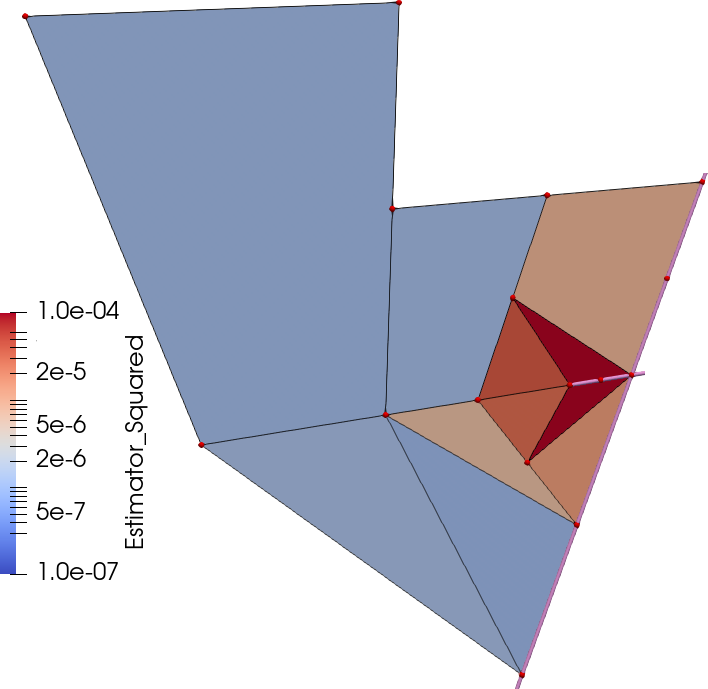}
        \caption{Case $1$ - \emph{New} algorithm}
        \label{fig:algo:aligned:example:new}
    \end{subfigure}\\
    \begin{subfigure}[!h]{0.32\textwidth}
        \centering
        \includegraphics[width=\textwidth]{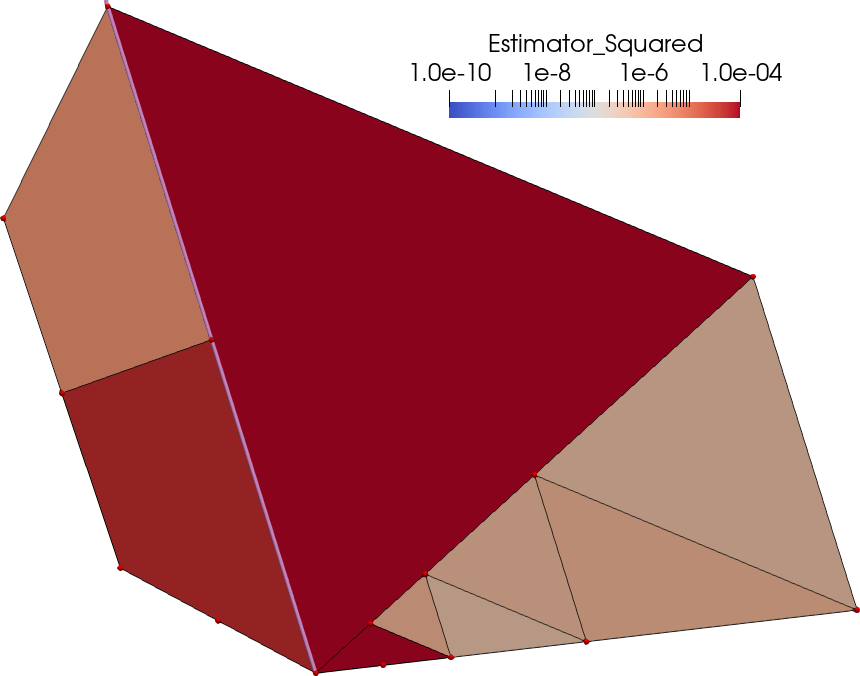}
        \caption{Case $2$ - \emph{Old} algorithm}
        \label{fig:algo:neigh:example:old}
    \end{subfigure}
    \begin{subfigure}[!h]{0.32\textwidth}
        \centering
        \includegraphics[width=\textwidth]{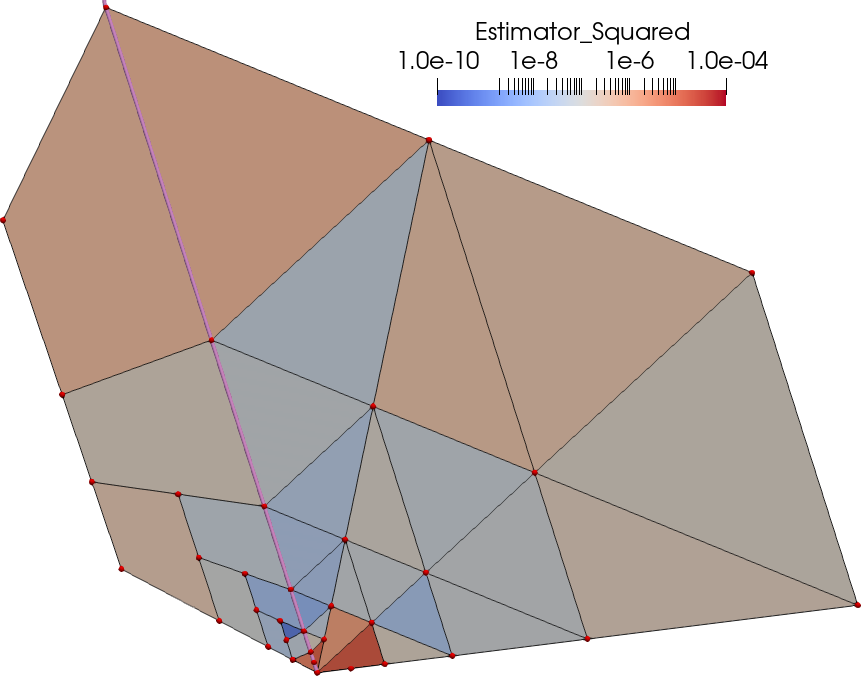}
        \caption{Case $2$ - \emph{New} algorithm}
        \label{fig:algo:neigh:example:new}
    \end{subfigure}
    \caption{Refinement of $E$ - Measure of the local error estimator on a numerical case}
    \label{fig:algo:example}
\end{figure}

The proposed \texttt{REFINE} algorithm represents an extension and enhancement of the concepts introduced in \cite{BERRONE2022103770, doi:10.1142/S0218202517500233}.
In the subsequent, we elucidate our motivation for developing Algorithm~\ref{alg:ref} by drawing a comparison with Algorithm~$4$ presented in \cite{BERRONE2022103770}.
In what follows, we reserve the label \emph{New} to the former algorithm, and we refer as the \emph{Old} to the latter one.

A pivotal innovation in our work is the incorporation of Check~\eqref{check:two} in Line~\ref{alg:ref:al} within the \texttt{CHECK-QUALITY} Algorithm~\ref{alg:quality}.
This subtle modification change is instrumental in controlling Condition~\eqref{eq:minimal:conditions:al}, resulting in a comprehensive enhancement of the mesh quality.
In Figures~\ref{fig:algo:aligned}, we provide an illustrative example where Check~\eqref{check:two} plays a fundamental role.

Consider a scenario where the triangle with aligned edges $E$ of Figure~\ref{fig:algo:aligned:original} is split in the edge $e$.
We assume $|\mathcal{I}_e| = \ell$, $\# \mathcal{I}_e = 3$, and the existence of a short edge in $\mathcal{I}_e$ of length $(1-p)\ell$, with $c_{\rho} = 0.5$.
In the \emph{Old} algorithm, when the portion $(1-p)$ is sufficiently small, i.e. $p \in (\frac{4}{5}, 1)$, Check~\eqref{check:one} permits the creation of a new aligned edge, as depicted in Figure~\ref{fig:algo:aligned:old}.
However, this new edge deteriorates Condition~\eqref{eq:minimal:conditions:al}.
On the other hand, in the \emph{New} algorithm, the introduction of Check~\eqref{check:two} enables the collapse of the cut direction to an original vertex, as illustrated in Figure~\ref{fig:algo:aligned:new}.
Figures~\ref{fig:algo:aligned:example:old}-\ref{fig:algo:aligned:example:new} show the local error estimators in a numerical example where this specific condition arises.
It is evident from the colouring of the cells that the \emph{Old} algorithm (Figure~\ref{fig:algo:aligned:example:old}) results in a higher mean error estimator compared to the \emph{New} proposed version, (Figure~\ref{fig:algo:aligned:example:new}) at the same refinement step $m$.

The \emph{extension} operation, a key component of the \texttt{REFINE} Algorithm~\ref{alg:ref}, is encapsulated in Line~\ref{alg:ref:extend} and represents the other noteworthy innovation in our work.
This operation addresses a crucial aspect of the refinement process by extending it to all the polygonal cells $E_q \in \mathcal{Q}$ containing marked edges that fail Checks~\eqref{check:one}-\eqref{check:two} in the \texttt{CHECK-QUALITY} algorithm.
Note that, during this \texttt{CHECK-QUALITY} task, we set $s=1$ since no edge split is necessary.
We recall that we selectively mark the edges $e \in E$ split by the direction $\mathbf{t}$.
This strategic approach allows us to handle the scenario depicted in Figures~\ref{fig:algo:neigh}.
We illustrate a situation where the estimator consistently requires the refinement of cell $E$ of Figure~\ref{fig:algo:neigh:original} in each step $m$.
The initial cut applies along the longest edge $e$.
Moreover, we suppose the neighbour cell $E_n$ to have $r_{E_n}$ larger than $e$ and $h_{E_n} = e$.
As the refined cells $E_1, E_2, \dots$ are always triangles, in the \emph{Old} algorithm proposed in \cite{BERRONE2022103770}, as seen in Figure~\ref{fig:algo:neigh:old}, the neighbour cell $E_n$ experiences the creation of a non-uniform aligned edge collection $\mathcal{I}^e_{E_n}$ as a consequence.
In contrast, our proposed \emph{New} version, shown in Figure~\ref{fig:algo:neigh:new}, extends the refinement if the quality Check~\ref{check:two} is not met.
In Figures~\ref{fig:algo:neigh:example:old}-\ref{fig:algo:neigh:example:new}, we present a numerical test to illustrate the local mean estimator substantial improvement in the \emph{New} case.
This empirical evidence reinforces the efficacy of our approach.

We conclude with a final consideration regarding the pivotal role played by the constants $c_{\rho}$ and $c_{al}$.
As extensively discussed in \cite{BERRONE2022103770}, the configuration of the resulting cells in the mesh $\mathcal{T}_{\Omega}^{m+1}$ is significantly influenced by the values assigned to $c_{\rho}$.
Indeed, the number of aligned edges experiences a notable increase when $c_{\rho}$ tends towards $0$.
Conversely, when $c_{\rho}$ assumes a much larger value, the prevalence of cells shape shifts towards triangular mesh elements.
The parameter $c_{al}$ plays a analogous role, with an increase in aligned edges as $c_{al}$ approaches zero.

The interdependence of $c_{\rho}$ and $c_{al}$ is evident, stemming from their inherent connection to the geometric properties of the cells.
When $c_{\rho} \gg 1.0$ and the mesh elements exhibit uniform sizes, Check~\eqref{check:one}   predominantly fails, rendering Check~\eqref{check:two} unused and rendering the value of $c_{al}$ irrelevant.
The converse holds true as well.
We are going to investigate and discuss the impact of $c_{\rho}$ and $c_{al}$ on the mesh element shapes in $\mathcal{T}_{\Omega}^m$ in the section devoted for numerical results.

\section{The problem and the discretization}
\label{sec:problem}
We are going to demonstrate the effectiveness of the newly proposed Process~\eqref{eq:schema} by focusing on the simulation of the flow in fractured porous media.
Specifically, we employ the Discrete Fracture Network (DFN) framework.
We select this problem due to its natural generation of an initial discretization $\mathcal{T}^0_{\Omega}$ that exhibits highly \emph{badly-shape} features. 
It is crucial to emphasize that the strategy outlined in this paper is applicable to any two-dimensional problem where polygonal meshes are valuable.

Consider a DFN domain $\Omega$, comprising the union of $I$ polygonal planar fractures $F_i \in \mathbb{R}^3$, $\Omega := \bigcup_i F_i$, with $i \in \{1, \dots, I\}$.
These fractures are randomly oriented in three-dimensional space, and their intersections create $M$ segments denoted by $S_m \in \mathbb{R}^3$, where $m \in \{1, \dots, M\}$.
We assume that each segment $S_m$ results from the intersection of exactly two fractures, establishing a one-to-one relationship between the intersection index $m$ and two fracture indices, i.e., $S_m = \bar{F_i} \cap \bar{F_j}$, $\forall m \in \{1, \dots, M\}$.

The boundary $\partial \Omega = \bigcup_i \partial F_i$ is divided into a Dirichlet portion $\Gamma_D \neq \emptyset \subseteq \partial \Omega$ and a Neumann portion $\Gamma_N = \partial \Omega \setminus \Gamma_D$.
Dirichlet conditions are prescribed by the function $g_D \in H^{\frac{1}{2}}(\Gamma_D)$ and a homogeneous Neumann condition is applied on $\Gamma_N$. 
We remark that different boundary conditions can be considered.

For each fracture $F_i$, we introduce the following functional spaces:
\begin{linenomath}
    \begin{equation*}
        V^D_i = \left\{ v \in H^1(F_i) : v_{|_{\Gamma_D \cap \partial F_i}} = g_D \right\}, \quad
        V_i = \left\{ v \in H^1(F_i) : v_{|_{\Gamma_D \cap \partial F_i}} = 0 \right\}.
    \end{equation*}
\end{linenomath}
Moreover, on the whole DFN $\Omega$ we define
\begin{linenomath}
    \begin{eqnarray}
        \label{eq:spaces}
        V^D &=& \left\{ v : v_{|_{F_i}} \in V^D_i\ \forall i \in \{1, \dots, I\},\ \gamma_{S_m}(v_{|_{F_i}}) = \gamma_{S_m}(v_{|_{F_j}})\ \forall m \in \{1, \dots, M\} \right\}, \\
        V &=& \left\{ v : v_{|_{F_i}} \in V_i\ \forall i \in \{1, \dots, I\},\ \gamma_{S_m}(v_{|_{F_i}}) = \gamma_{S_m}(v_{|_{F_j}})\ \forall m \in \{1, \dots, M\} \right\}\nonumber,
    \end{eqnarray}
\end{linenomath}
where $\gamma_{S_m} : H^1(F_i) \to H^{\frac{1}{2}}(S_m)$ represents the trace operator onto $S_m$, $\forall m \in \{1, \dots, M\}$.

We seek the resolution of the following symmetric PDE: find $U \in V^D$ such that:
\begin{linenomath}
    \begin{equation}
        \label{eq:prob:cont}
        \sum_{F_i} \left(K_i \nabla U_i, \nabla v \right)_{F_i} = \sum_{F_i} \left(Q_i, v\right)_{F_i} \ \forall v \in V.
    \end{equation}
\end{linenomath}
This equation models Darcy's law applied to the DFN.
In this context, $U$ represents the hydraulic head on the network, $U_i$ its restriction on $F_i$, $K_i \in \mathbb{R}$ denotes the fracture transmissivity, and $Q_i \in L^2 (F_i)$ acts as the fracture source term.

\subsection{DFN Discretization}
The application of the Process~\eqref{eq:schema} and the resolution to Problem~\eqref{eq:prob:cont} with the VEM requires the ability to construct a partition $\mathcal{T}^0_{\Omega}$ on the DFN that globally conforms to the fracture intersections.
Our approach is based on the methodology detailed in \cite{BENEDETTO201623}.
The generated mesh $\mathcal{T}^0_{\Omega}$ results from the following steps:
\begin{itemize}
    \item On each fracture $F_i$, we generate a local tessellation $\mathcal{T}^{loc}_{F_i}$.
    This discretization is not necessarily conforming to the fracture intersections and it is entirely independent of other DFN fractures.
    \item On each $F_i$, we cut the elements of $\mathcal{T}^{loc}_{F_i}$ using the fracture intersection segments $S_m$, $\forall S_m \subset F_i$. If a segment $S_m$ terminates inside a polygonal cell of $\mathcal{T}^{loc}_{F_i}$, the segment is extended to the end of the cell. It is important to note that this extension does not alter the original domain $\Omega$ shape. This process yields a new polygonal local mesh $\mathcal{T}^{cnf}_{F_i}$ that conforms locally to the fracture intersections.
    \item On each intersection $S_m = \bar{F_i} \cap \bar{F_j}$, $m \in \{1, \dots, M\}$, we unify all the mesh nodes of $T^{cnf}_{F_i}$ and $T^{cnf}_{F_j}$ lying on $S_m$. The resulting union points are added to the edges lying on $S_m$ of each mesh $T^{cnf}_{F_i}$ and $T^{cnf}_{F_j}$.
\end{itemize}
The union of the $T^{cnf}_{F_i}$, $\forall i \in \{1, \dots, I\}$ meshes produces the globally conforming mesh $\mathcal{T}^0_{\Omega}$.

In this work, we opt for the fracture \emph{minimal mesh} as our choice for $\mathcal{T}^{loc}_{F_i}$, as introduced in Section~\ref{sec:setting}.
This entails the use of a discretization formed by the convex polygonal approximation of $\bar{F_i}$. 
We will henceforth refer to $\mathcal{T}^0_{\Omega}$ as the \emph{minimal mesh} for the DFN case.

It is important to note that the mesh generation process described naturally translates to the generation of a tessellation with \emph{bad-shape} elements, featuring generic polygonal elements and aligned edges.
The quality of the resulting $\mathcal{T}^0_{\Omega}$ is strongly influenced by the fractures' spatial location in the tree-dimensional space.
In real DFNs, $\mathcal{T}^0_{\Omega}$ exhibits numerous aligned edges and neighbouring polygonal cells with significant variation in size.
We rely on the refining strategy proposed here to enhance the mesh quality.

Thanks to the global conformity of the mesh $\mathcal{T}^0_{\Omega}$ and the ability to identify the DOFs on each segment $S_m$, where $m \in \{1, \dots, M\}$, we can define the global discrete spaces of Equation~\eqref{eq:spaces} as
\begin{linenomath}
    \begin{eqnarray}
        \label{eq:spaces:discrete}
        V^D_k &=& \left\{ v \in V^D : v_{|_{F_i}} \in C^0(F_i),\ v_{|_E} \in V^E_k,\ \forall E \in F_i,\ \forall i \in \{1, \dots, I\}\right\}, \\
        V_k &=& \left\{ v \in V : v_{|_{F_i}} \in C^0(F_i),\ v_{|_E} \in V^E_k,\ \forall E \in F_i,\ \forall i \in \{1, \dots, I\}\right\}.\nonumber
    \end{eqnarray}
\end{linenomath}
We refer to Equation~\eqref{eq:space:discrete:local} for the definition of $V^E_k$ on each element $E \in \mathcal{T}_{\Omega}$.

We emphasize that all the quantities for $E \in \mathcal{T}_{\Omega}$ are computed with respect to the reference system tangential to the fracture $F_i$ to which the element $E$ belongs.
Despite the DFN being immersed in $\mathbb{R}^3$, the use of the reference tangential system allows us to work exclusively in the geometric dimension $2$.
Additionally, the reader shall notice that, in Equation~\eqref{eq:spaces:discrete}, we employ the same symbol $V_k$ of the space as in Equation~\eqref{eq:space:discrete}, as it naturally extends this to the DFN case.

We discretize Problem~\eqref{eq:prob:cont} as follows: find $u \in V^D_k$ such that: 
\begin{linenomath}
    \begin{equation}
        \label{eq:prob:disc}
        \sum_{F_i} \left(K_i \Pi^0_{k-1}\nabla u_i, \nabla \Pi^0_{k-1} v \right)_{F_i} + S^E_{F_i}(u_i - \Pi^{\nabla}_k u_i, v - \Pi^{\nabla}_k v) = \sum_{F_i} \left(Q_i, \Pi^0_{k-1} v\right)_{F_i} \ \forall v \in V_k
    \end{equation}
\end{linenomath}
where $S^E_{F_i}$ is the stabilizing bilinear form. 
Among all the available choices, see~\cite{doi:10.1142/S021820251750052X}, we opt for the standard form:
\begin{linenomath}
    \begin{equation*}
        S^E_{F_i}(u, v) = K_i \sum_{j=1}^{\# V^E_k} \text{dof}^E_j(u) \text{dof}^E_j(v)\quad \forall u \in V^D_k, v \in V_k.
    \end{equation*} 
\end{linenomath}

\subsection{DFN Mesh Adaptivity}
We introduce the VEM residual-type error estimator for the DFN case, which we use in the \texttt{ESTIMATE} function of the adaptive Process~\eqref{eq:schema}.
We follow the analysis proposed in \cite{doi:10.1142/S0218202517500233} and its extension to the DFN case reported in \cite{BERRONE2021103502}.

Given $u^m$ as the discrete solution generated by the \texttt{SOLVE} task on each iteration $m$ of Process~\eqref{eq:schema}, we define the function $\eta_{\Omega} : V^D_k \to \mathbb{R}$ s.t.
\begin{linenomath}
    \begin{equation}
        \label{eq:est}
        \eta_{\Omega}^2(u^m) :=  \sum_{E \in \mathcal{T}^m_{\Omega}} \frac{D^2_E}{K_i}\left|\left|\Pi^0_{k-1} Q_i + K_i \Delta u^{\pi}_i\right|\right|_E^2 + \sum_{e \in \mathcal{E}^m_{\Omega} \setminus \Gamma_D} \frac{|e|}{K_e} \left|\left| \sum_{\tilde{E} \in \mathcal{N}_e} K_i \nabla u^{\pi}_i \cdot n^{\tilde{E}}_e \right|\right|_e^2 + \sum_{E \in \mathcal{T}^m_{\Omega}} \frac{D^2_E}{K_i}\left|\left|Q_i - \Pi^0_{k-1} Q_i\right|\right|_E^2,
    \end{equation}
\end{linenomath}
where $u^{\pi}_i := (\Pi^{\nabla}_k u^m)_{|_{F_i}}$, $n^{\tilde{E}}_e$ is the unit normal vector to $e$ pointing outward with respect to $\tilde{E} \in \mathcal{N}_e$ and $K_e := \sum_{\tilde{E} \in \mathcal{N}_e} K_i$.
The subscript $i$ refers to the fracture $F_i$ to which the element $E$ belongs.

It is possible to prove, as discussed in \cite{doi:10.1142/S0218202517500233}, that there exist two positive constants $C_*$ and $C^*$ such that:
\begin{linenomath}
    \begin{equation*}
        C_* \cdot \eta_{\Omega}(u^m) \leq \left|\left|\left| U - \Pi^{\nabla}_k u^m \right|\right|\right|_{\Omega} \leq C^* \cdot \eta_{\Omega}(u^m),
    \end{equation*}
\end{linenomath}
where $U$ is the solution to Problem~\eqref{eq:prob:cont}, and $u^m$ is the solution to Problem~\eqref{eq:prob:disc}.
For error control, we use the energy norm $|||\cdot||| : V^D \to \mathbb{R}$,
\begin{linenomath}
    \begin{equation*}
        \left|\left|\left| v \right|\right|\right|_{\Omega} := \sup_{w \in V} \frac{\sum_{F_i} \left(K_i \nabla v, \nabla w \right)_{F_i}}{\left(\sum_{F_i} \left|\left|\sqrt{K_i} \nabla w \right|\right|^2_{F_i}\right)^{\frac{1}{2}}} \quad \forall v \in V^D.
    \end{equation*}
\end{linenomath}

Finally, on each polygon $E \in \mathcal{T}^m_{\Omega}$, we introduce the local estimator $\eta_E : V^D_k \to \mathbb{R}$ such that:
\begin{linenomath}
    \begin{equation}
        \label{eq:est:local}
        \eta_{E}^2(u^m) := \frac{D^2_E}{K_i}\left|\left|\Pi^0_{k-1} Q_i + K_i \Delta u^{\pi}_i\right|\right|_E^2 + \sum_{e \in \mathcal{E}_E \setminus \Gamma_D} \frac{|e|}{\# \mathcal{N}_e \cdot K_e} \left|\left| \sum_{\tilde{E} \in \mathcal{N}_e} K_i \nabla u^{\pi}_i \cdot n^{\tilde{E}}_e \right|\right|_e^2 + \frac{D^2_E}{K_i}\left|\left|Q_i - \Pi^0_{k-1} Q_i\right|\right|_E^2.
    \end{equation}
\end{linenomath}

\begin{remark}[Local flux estimation]
    We incorporate the inverse of $\# \mathcal{N}_e$ to weight the local flux estimation. This choice mirrors the constant $\frac{1}{2}$ typically used in classical $\mathbb{R}^2$ applications, given that $\# \mathcal{N}_e = 2$. 
\end{remark}
\section{Numerical Results}
\label{sec:results}
In this section, we analyse the performances of the proposed adaptive Scheme~\eqref{eq:schema}.
For all the performed tests, the initial mesh $\mathcal{T}_{\Omega}^0$ is selected as the \emph{minimal mesh} introduced in the previous sections.

The marking strategy employed is based on the approach proposed in \cite{doi:10.1137/0733054}.
The parameter $\theta \in [0.0, 1.0]$ is used to identify a subset of cells $\mathcal{M}^m \subseteq \mathcal{T}_{\Omega}^m$ such that:
\begin{linenomath}
    \begin{equation*}
        \sum_{E \in \mathcal{M}^m} \eta_E^2(u^m) \geq \theta \eta_{\Omega}^2(u^m),
    \end{equation*}
\end{linenomath}
where $\eta_{\Omega}$ and $\eta_E$ are the estimators introduced in Equation~\eqref{eq:est} and Equation~\eqref{eq:est:local}, respectively.
Throughout all the tests, we set $\theta = 0.5$, as it yields the most favourable results.

To gauge the mesh $\mathcal{T} = \mathcal{T}^m_{\Omega}$ quality performance, we introduce the following indicators:
\begin{itemize}
    \item $AR^{Rr}_E := \frac{R_E}{r_E}$ and $AR^{Hr}_E := \frac{R_E}{h_E}$, serving as local approximations of the inverses of constants $\gamma_r$ and $\gamma_h$, respectively;
    \item $\# \Delta_{\mathcal{T}}$ and $\# \diamond_{\mathcal{T}}$ represent the number of cells $E \in \mathcal{T}$ with $3$ and $4$ vertices, respectively. These values measure the presence of ``real'' triangles and quadrilateral on the tessellation.
    \item $\# \Delta^{al}_{\mathcal{T}}$ and $\# \diamond^{al}_{\mathcal{T}}$ denote the number of cells $E \in \mathcal{T}$ whose polygons $\hat{E}$, obtained by unifying the aligned edges $e \in \mathcal{E}_E$, are formed by $3$ and $4$ vertices, respectively. Note that, $\Delta^{al}_{\mathcal{T}}$ includes also $\Delta_{\mathcal{T}}$ cells and some cells of $\diamond_{\mathcal{T}}$; similarly, $\diamond^{al}_{\mathcal{T}}$ contains the resulting part of the cells from $ \diamond_{\mathcal{T}}$.
    \item $R^{\Delta} := \frac{\# \Delta_{\mathcal{T}}}{\# \mathcal{T}}$, $R^{\diamond} := \frac{\# \diamond_{\mathcal{T}}}{\# \mathcal{T}}$, $R^{P} := 1 - (R^{\Delta} + R^{\diamond})$ measure the fraction of different polygons contained in the tessellation.
    \item $R_{al}^{\Delta} := \frac{\# \Delta^{al}_{\mathcal{T}}}{\# \mathcal{T}}$, $R_{al}^{\diamond} := \frac{\# \diamond^{al}_{\mathcal{T}}}{\# \mathcal{T}}$, measure the percentage of different shapes contained in the tessellation, accounting only polygons with not aligned edges.
    \item $\text{Ef }_{\# \mathcal{T}}^{\text{DOFs}} := \frac{\# \mathcal{T}}{\# DOFs}$ computes the efficiency of the mesh elements in approximating the numerical solution. A large value indicates better efficiency.
\end{itemize}

\subsection{Test 1: Classic bi-dimensional Domain}
\begin{figure}[!h]
    \centering
    \begin{subfigure}[]{0.32\textwidth}
        \centering
        \includegraphics[width=\textwidth]{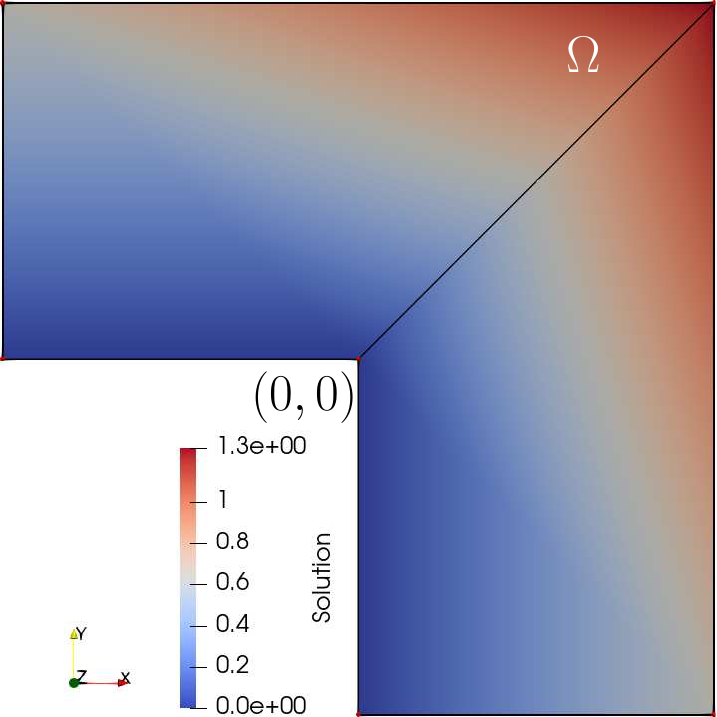}
        \caption{$\mathcal{T}^0_{\Omega}$}
        \label{fig:lshape:network:t0:05}
    \end{subfigure}
    \begin{subfigure}[]{0.32\textwidth}
        \centering
        \includegraphics[width=\textwidth]{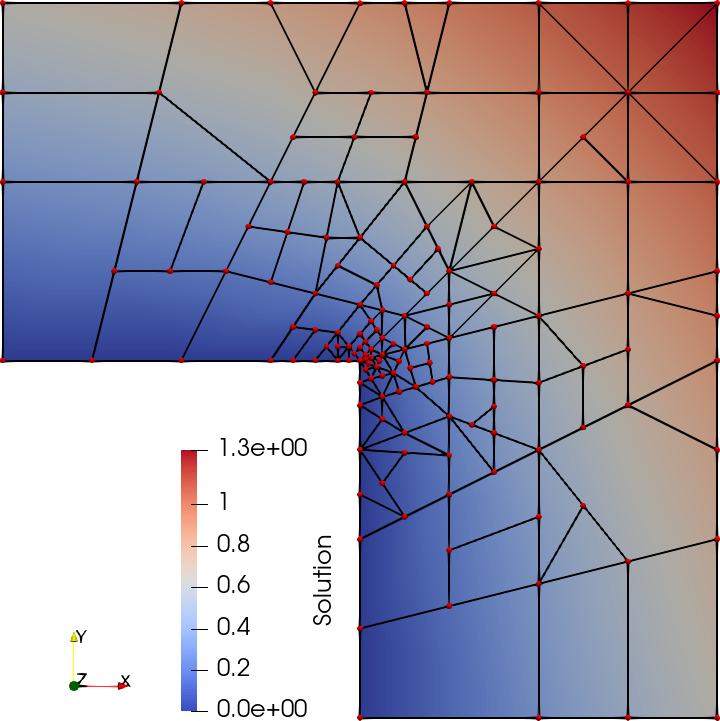}
        \caption{Step $m=15$, $c_\rho = 0.5$}
        \label{fig:lshape:network:t15:05}
    \end{subfigure}
    \begin{subfigure}[]{0.32\textwidth}
        \centering
        \includegraphics[width=\textwidth]{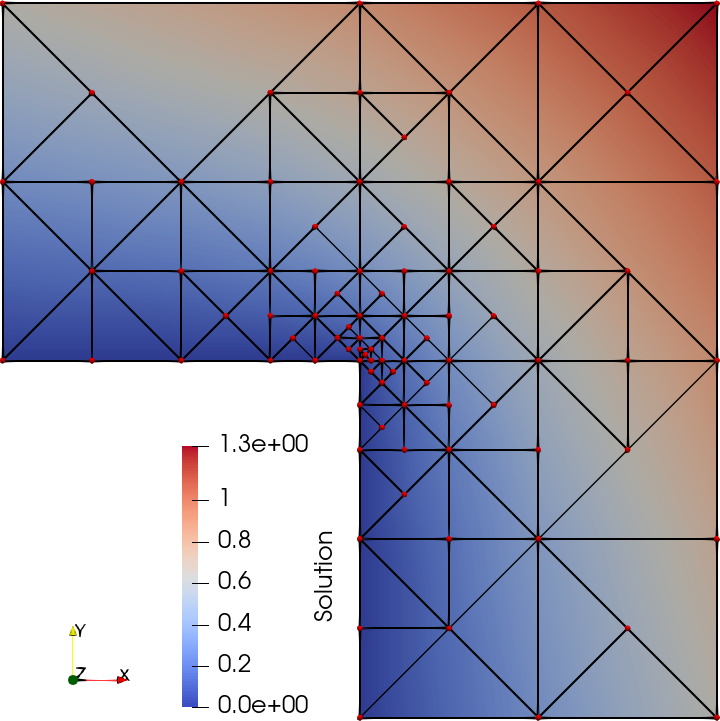}
        \caption{Step $m=15$, $c_\rho = 1.5$}
        \label{fig:lshape:network:t15:15}
    \end{subfigure}
    \caption{Test $1$ - Domain with solution and, $k = 1$, $c_{al} = 1.0$.}
    \label{fig:lshape:network}
\end{figure}
\begin{figure}[!h]
    \centering
    \begin{subfigure}[]{0.32\textwidth}
        \centering
        \includegraphics[width=\textwidth]{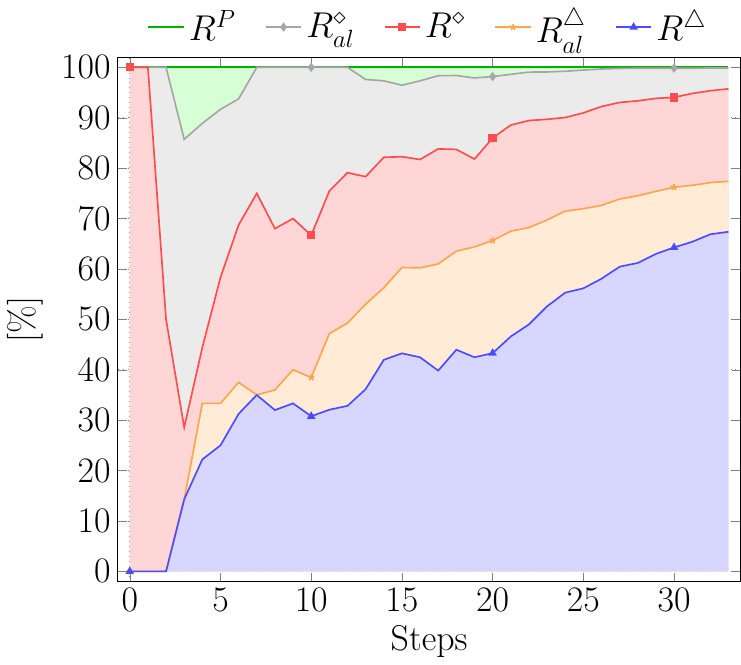}
        \caption{$R^{\star}$ vs Steps, $c_\rho = 0.5$}
        \label{fig:lshape:m_quality:tri:05}
    \end{subfigure}
    \begin{subfigure}[]{0.32\textwidth}
        \centering
        \includegraphics[width=\textwidth]{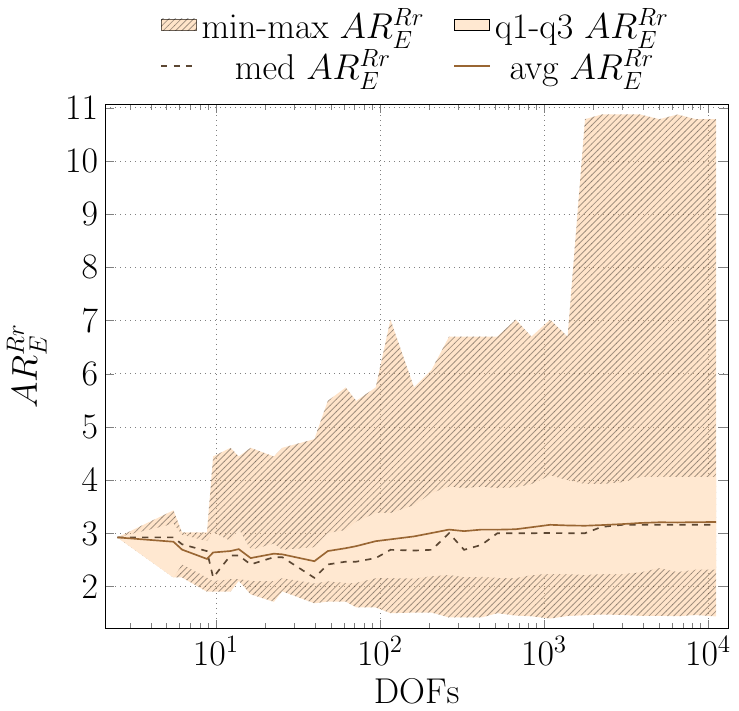}
        \caption{$AR^{Rr}_E$ vs DOFs, $c_{\rho} = 0.5$}
        \label{fig:lshape:dofs_ar:Rr:05}
    \end{subfigure}
    \begin{subfigure}[]{0.32\textwidth}
        \centering
        \includegraphics[width=\textwidth]{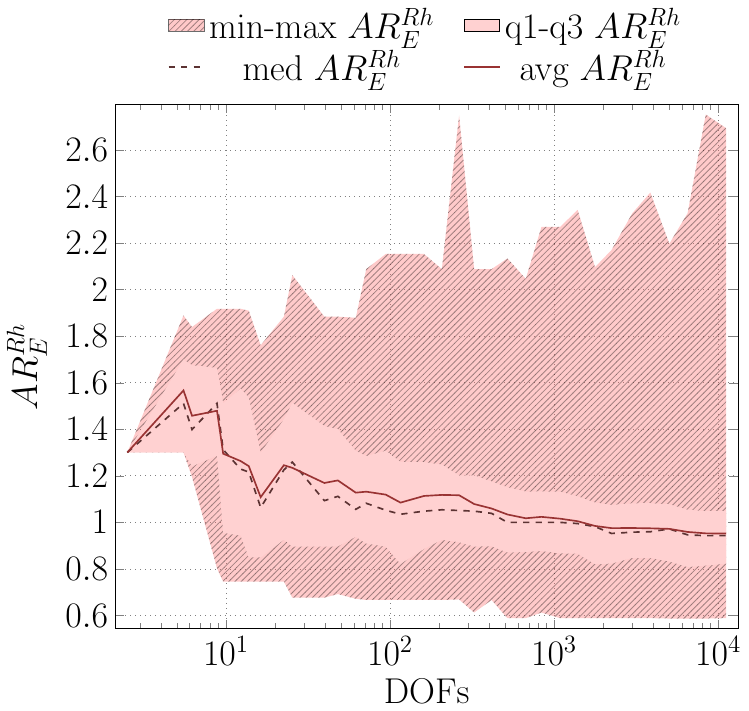}
        \caption{$AR^{Rh}_E$ vs DOFs, $c_{\rho} = 0.5$}
        \label{fig:lshape:dofs_ar:Rh:05}
    \end{subfigure}
    \begin{subfigure}[]{0.32\textwidth}
        \centering
        \includegraphics[width=\textwidth]{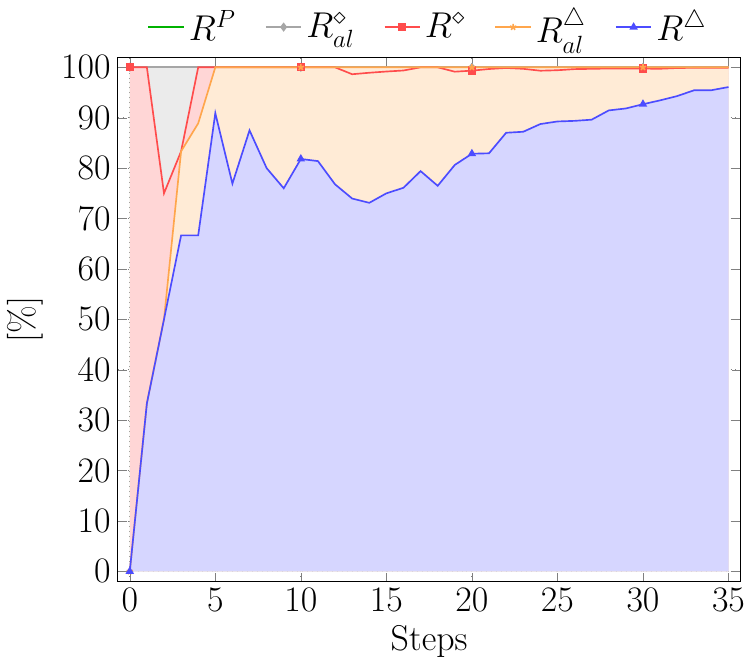}
        \caption{$R^{\star}$ vs Steps, $c_\rho = 1.5$}
        \label{fig:lshape:m_quality:tri:15}
    \end{subfigure}
    \begin{subfigure}[]{0.32\textwidth}
        \centering
        \includegraphics[width=\textwidth]{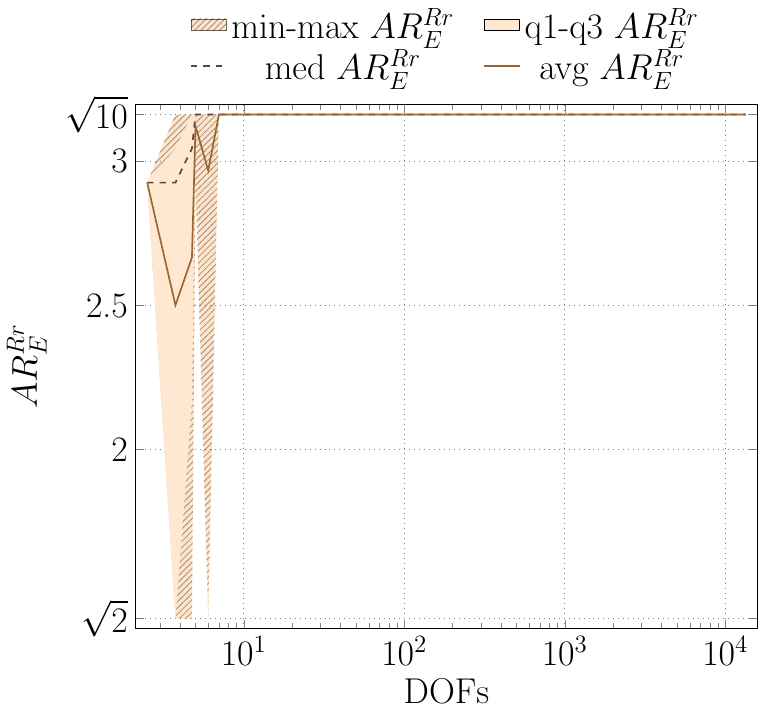}
        \caption{$AR^{Rr}_E$ vs DOFs, $c_{\rho} = 1.5$}
        \label{fig:lshape:dofs_ar:Rr:15}
    \end{subfigure}
    \begin{subfigure}[]{0.32\textwidth}
        \centering
        \includegraphics[width=\textwidth]{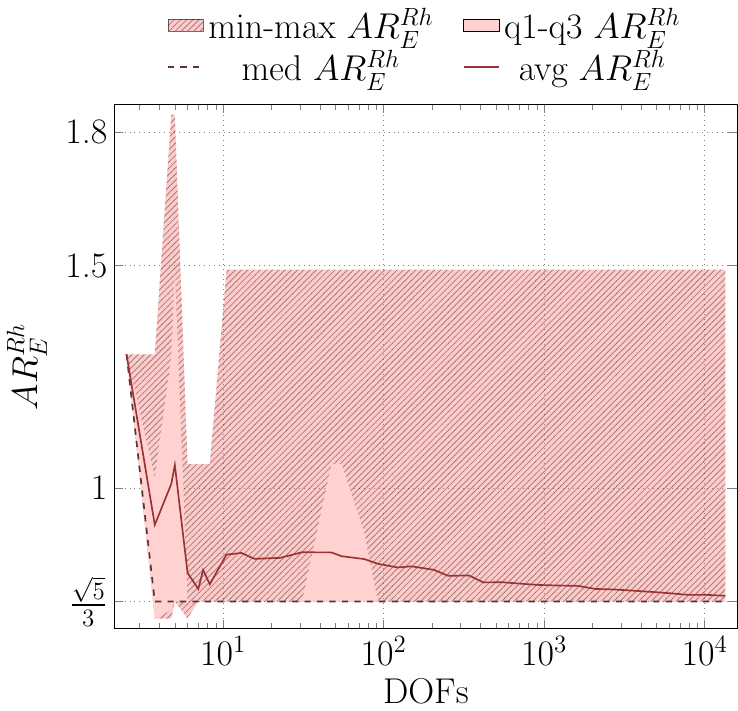}
        \caption{$AR^{Rh}_E$ vs DOFs, $c_{\rho} = 1.5$}
        \label{fig:lshape:dofs_ar:Rh:15}
    \end{subfigure}
    \caption{Test $1$ - mesh statistics, $k=1$, $c_{al} = 1.0$.}
    \label{fig:lshape:stat:15}
\end{figure}
\begin{figure}
    \centering
    \includegraphics[width=0.28\textwidth]{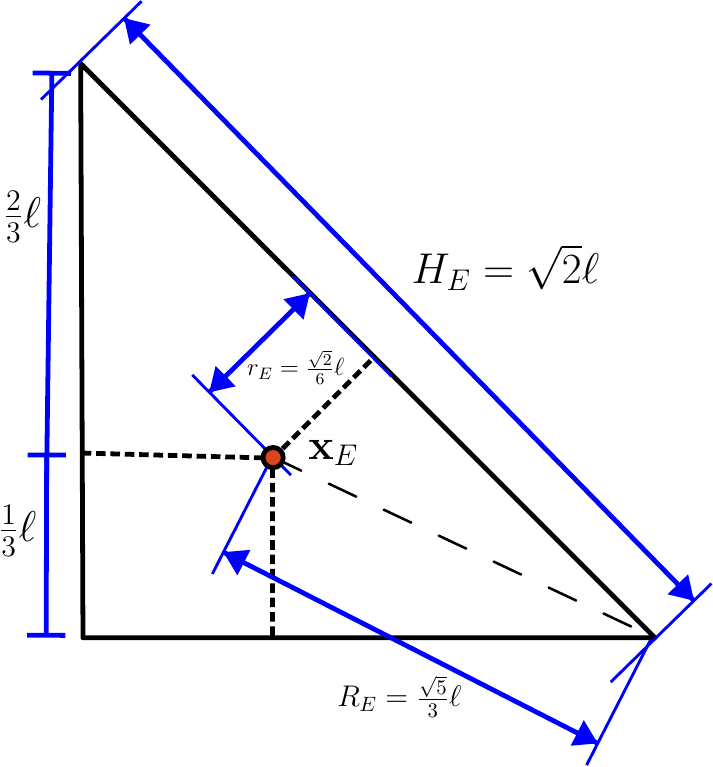}
    \caption{Ref. Triangle: $AR^{Rr}_E = \sqrt{10}$, $AR^{Rh}_E = \sqrt{5}/3$.}
    \label{fig:lshape:reference}
\end{figure}
\begin{figure}[!h]
    \centering
    \begin{subfigure}[]{0.32\textwidth}
        \centering
        \includegraphics[width=\textwidth]{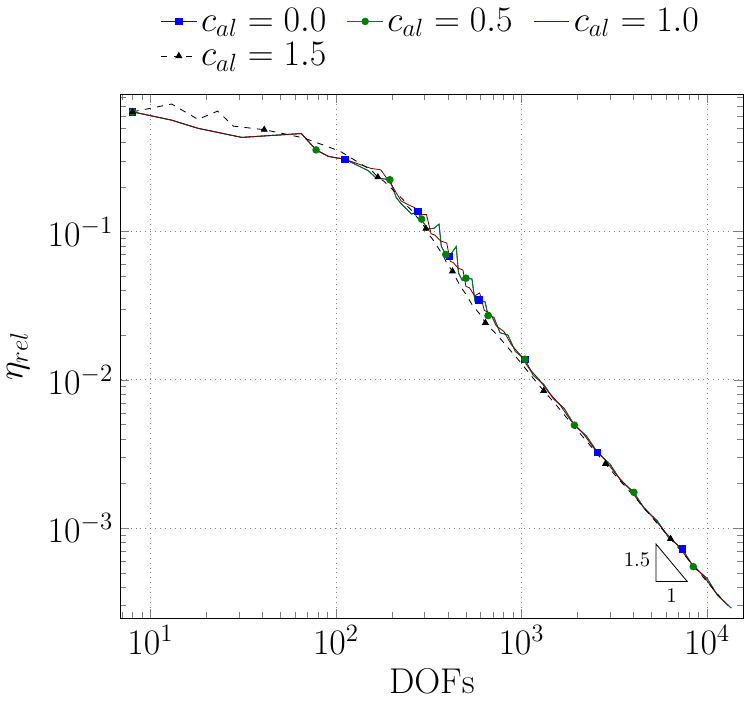}
        \caption{Estimator vs DOFs, $k=3$}
        \label{fig:lshape:est:al:dof:15}
    \end{subfigure}
    \begin{subfigure}[]{0.32\textwidth}
        \centering
        \includegraphics[width=\textwidth]{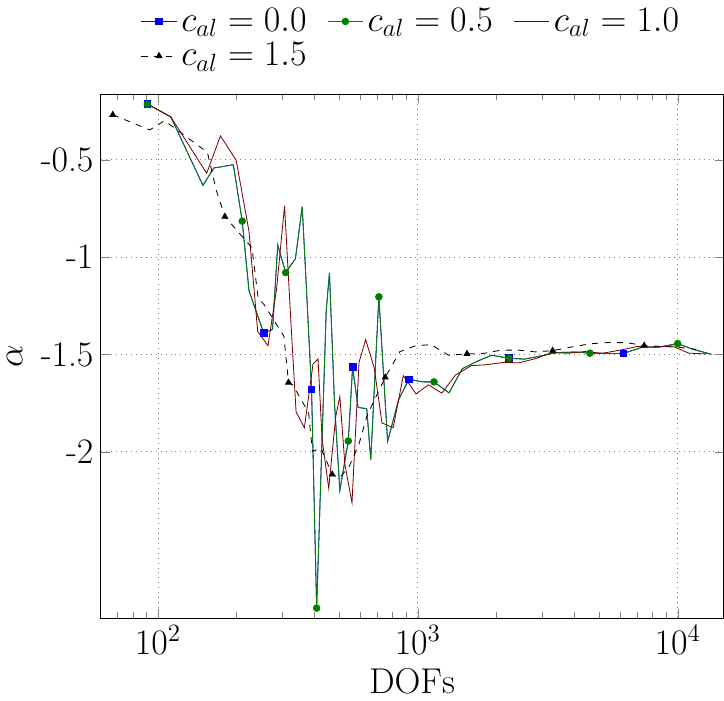}
        \caption{Slope $\alpha$ vs DOFs, $k=3$}
        \label{fig:lshape:est:al:slope:15}
    \end{subfigure}
    \begin{subfigure}[]{0.32\textwidth}
        \centering
        \includegraphics[width=\textwidth]{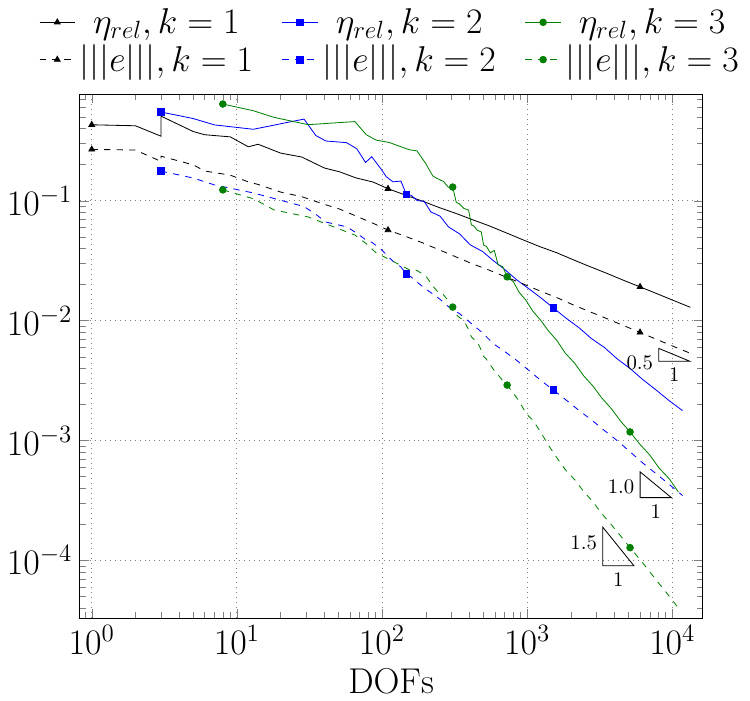}
        \caption{Estimator and Error vs DOFs, $c_{al} = 1.0$}
        \label{fig:lshape:est:dof_est:15}
    \end{subfigure}
    \caption{Test $1$ - convergence analysis, $c_{\rho} = 1.5$.}
    \label{fig:lshape:est:al:15}
\end{figure}
\begin{figure}[!h]
    \centering
    \begin{subfigure}[]{0.32\textwidth}
        \centering
        \includegraphics[width=\textwidth]{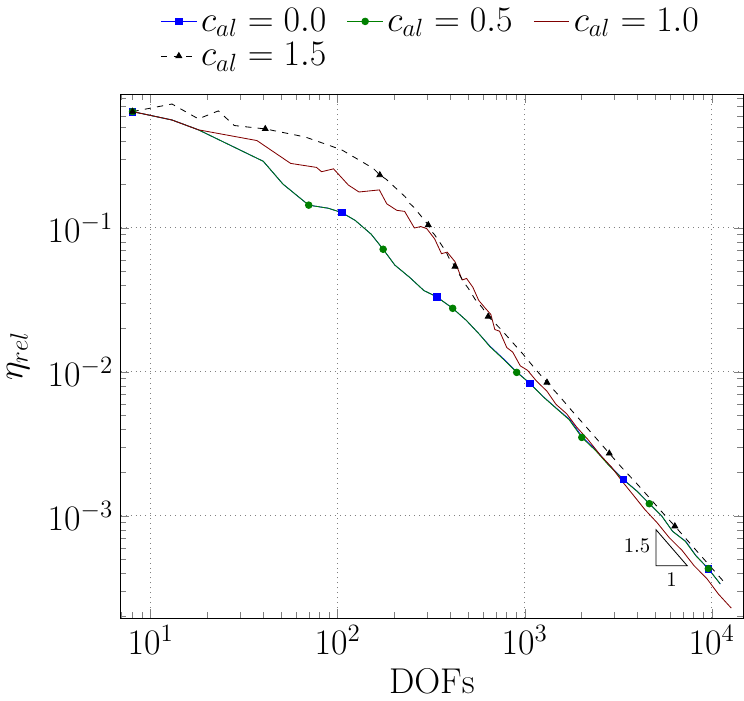}
        \caption{Estimator vs DOFs}
        \label{fig:lshape:est:al:dof:05}
    \end{subfigure}
    \begin{subfigure}[]{0.32\textwidth}
        \centering
        \includegraphics[width=\textwidth]{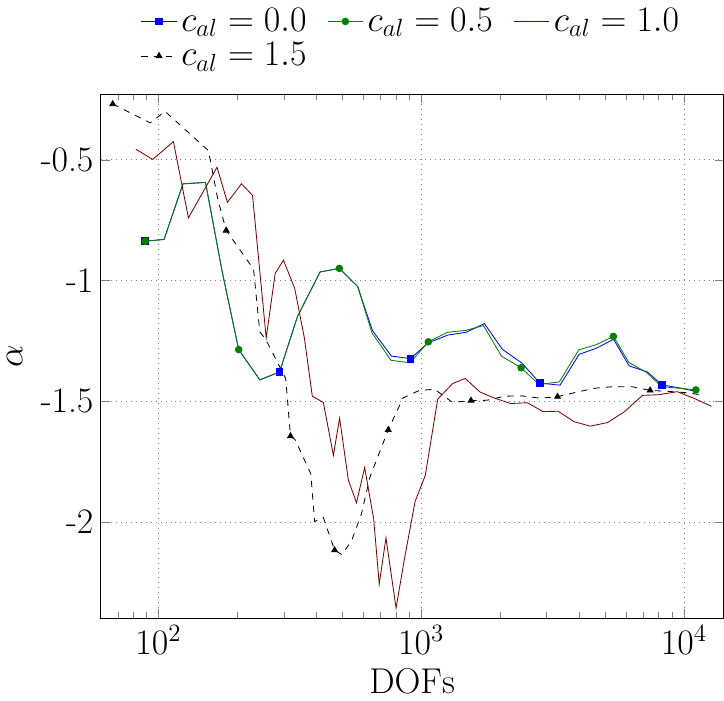}
        \caption{Slope $\alpha$ vs DOFs}
        \label{fig:lshape:est:al:slope:05}
    \end{subfigure}
    \begin{subfigure}[]{0.32\textwidth}
        \centering
        \includegraphics[width=\textwidth]{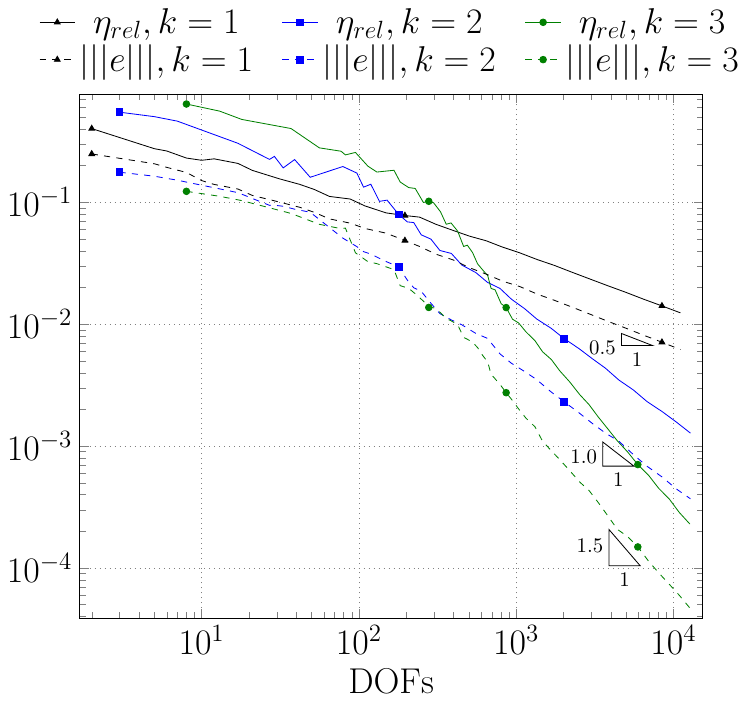}
        \caption{Estimator and Error vs DOFs, $c_{al} = 1.0$}
        \label{fig:lshape:est:dof_est:05}
    \end{subfigure}
    \begin{subfigure}[]{0.36\textwidth}
        \centering
        \includegraphics[width=\textwidth]{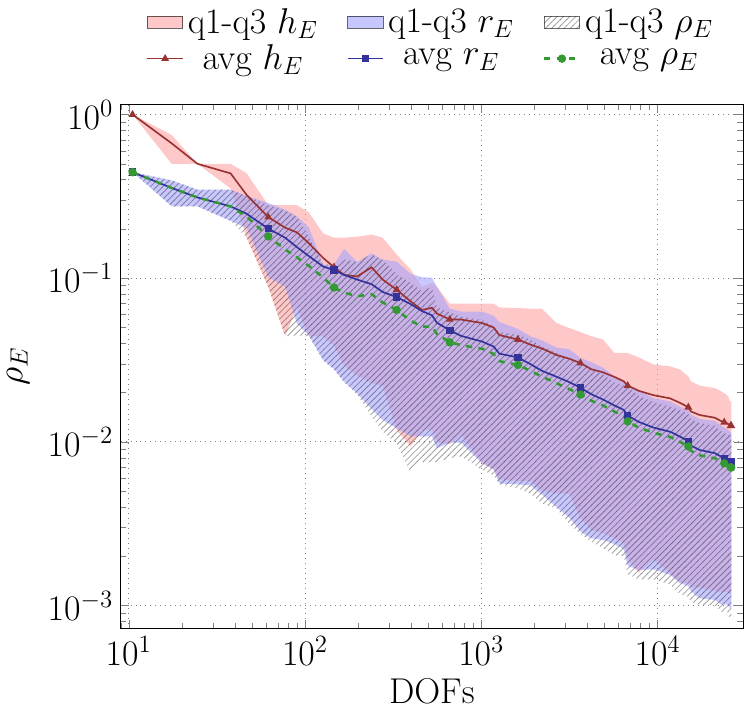}
        \caption{$\rho_E$ vs DOFs, $c_{al} = 0.5$}
        \label{fig:lshape:est:al:qE:05}
    \end{subfigure}
    \begin{subfigure}[]{0.36\textwidth}
        \centering
        \includegraphics[width=\textwidth]{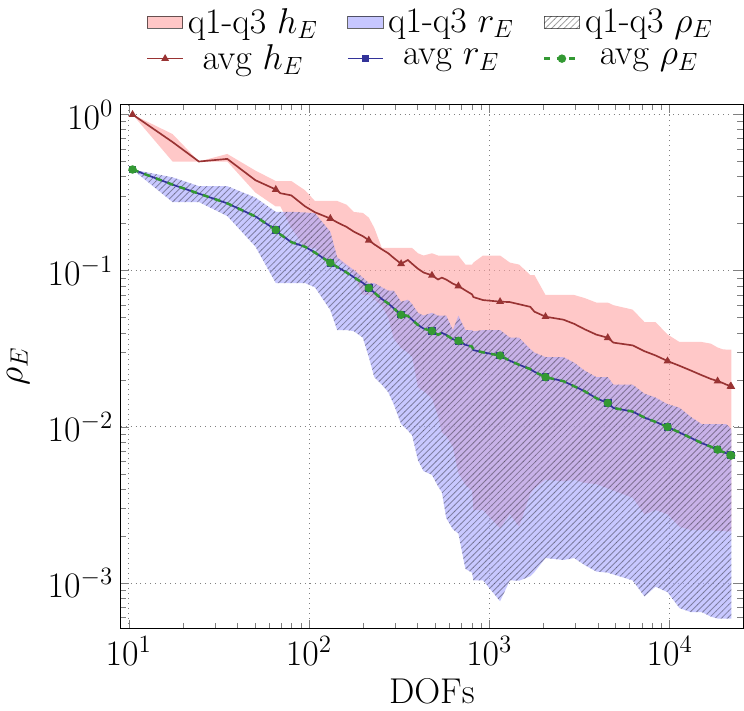}
        \caption{$\rho_E$ vs DOFs, $c_{al} = 1.0$}
        \label{fig:lshape:est:al:qE:10}
    \end{subfigure}
    \caption{Test $1$ - $c_{al}$ analysis, $c_{\rho} = 0.5$, $k=3$.}
    \label{fig:lshape:est:al:05}
\end{figure}
As first example, we present the solution to Problem~\eqref{eq:prob:disc} with $I=1$ on an L-shape domain $\Omega := (-1,1) \setminus (-1, 0)$.
This classic test serves as a measure of the application of the newly proposed refinement scheme to a generic two-dimensional problem.
The exact solution, prescribed with Dirichlet conditions, is given by $H(r, \beta) = r^{\frac{2}{3}}\sin{\frac{2}{3} (\beta + \frac{\pi}{2})}$, where $r$ and $\beta$ denote the polar coordinates.
We recall that the exact solution exhibits a singularity at the origin of the axis, as detailed in \cite{doi:10.1137/1.9781611972030.ch4}.
Figure~\ref{fig:lshape:network:t0:05} displays the domain with the initial \emph{minimal mesh} $\mathcal{T}^0_{\Omega}$.

We analyse the mesh quality and the solution approximation obtained during the refinement process for VEM order $k \in \{1, 2, 3\}$, $c_{\rho} \in \{0.5, 1.5\}$, ( as in \cite{BERRONE2022103770}), and $c_{al} \in \{0.0, 0.5, 1.0, 1.5\}$.
The adaptive procedure~\eqref{eq:schema} is interrupted when the total number of DOFs reaches $10^4$, irrespectively of the value of a-posteriori error $\eta_{\Omega}$.

Figure~\ref{fig:lshape:network:t15:05} and Figure~\ref{fig:lshape:network:t15:15} illustrate the mesh $\mathcal{T}_{\Omega}^m$ after $m=15$ steps for $c_{\rho} = 0.5$ and $c_{\rho} = 1.5$, respectively.
We report only the case with $c_{al} = 1.0$ and $k=1$, as other values exhibits similar attitude.
As emphasized in Section~\ref{sec:ref}, a value of $c_{\rho} < 1$ results in mesh elements with heterogeneous shapes, while $c_{\rho} > 1$ leads to a more shape-uniform mesh with a higher prevalence of triangular resulting cells.
The plots of Figure~\ref{fig:lshape:m_quality:tri:05} and Figure~\ref{fig:lshape:m_quality:tri:15} confirm these observations.
In Plot~\ref{fig:lshape:m_quality:tri:05}, we observe a combination of triangles and quadrilaterals ($R^{\Delta} + R^{\diamond}$ $\approx$ $60 \%$) and a high percentage of elements with aligned edges ($R^{\Delta}_{al} + R^{\diamond}_{al}$ $\approx$ $40 \%$).
On the contrary, Plot~\ref{fig:lshape:m_quality:tri:15} illustrate that the majority of $\mathcal{T}^m_{\Omega}$ cells are triangles ($R^{\Delta} \approx 70 \%$) or nearly triangles ($R^{\Delta}_{al} \approx 30 \%$) after approximately five steps.

In Figure~\ref{fig:lshape:dofs_ar:Rr:05}, Figure~\ref{fig:lshape:dofs_ar:Rr:15} and Figure~\ref{fig:lshape:dofs_ar:Rh:05}, Figure~\ref{fig:lshape:dofs_ar:Rh:15} we conduct a statistical analysis of the aspect ratios $AR^{Rr}_E$ and $AR^{Rh}_E$, respectively, on the mesh elements $E \in \mathcal{T}^m_{\Omega}$ for both the values $c_{\rho} = 0.5$ and $c_{\rho} = 1.5$.
For each iteration $m$, we present the minimum and maximum values (min-max $AR_E^{\star}$ area), the first and third quartile values (q1-q3 $AR_E^{\star}$ area), the median value (med $AR_E^{\star}$ curve), and the average value (avg $AR_E^{\star}$ curve) of $AR_E^{\star}$, $\star \in \{Rr, Rh\}$.
The statistical data reveals that the mesh $\mathcal{T}_{\Omega}^m$ obtained with $c_{\rho} = 0.5$ exhibits stable quality as the deviation from the mean value remains constant after a few iterations for both the quality $AR^{\star}_E$ indicators.
In contrast, elements of the mesh $\mathcal{T}_{\Omega}^m$ obtained with $c_{\rho} = 1.5$ resemble a rescaled version of the reference triangle in Figure~\ref{fig:lshape:reference}.
Indeed, the statistics converge to the characteristic quantities $AR_E^{Rr} = \sqrt{10}$ and $AR_E^{Rh} = \frac{\sqrt{5}}{3}$ of the reference triangle.

As observed in Section~\ref{sec:ref}, when $c_{\rho} > 1.0$ and the mesh elements present uniform sizes, the parameter $c_{al}$ becomes less relevant.
Indeed, in the \texttt{CHECK-QUALITY} Algorithm~\ref{alg:quality}, the Check~\eqref{check:one} often fails, and the Check~\eqref{check:two} is rarely used.
Figures~\ref{fig:lshape:est:al:15} confirm this observation, presenting an analysis of the parameter $c_{al} \in \{0.0, 0.5, 1.0, 1.5\}$ with a fixed value of $c_{\rho} = 1.5$.
Values $c_{al} = 0.0$ and $c_{al} = 1.5$ represent extremal cases.
Value $c_{al} = 0.0$ mimics the scenario where the Check~\eqref{check:two} of \texttt{CHECK-QUALITY} Algorithm~\ref{alg:quality} is disabled.
On the other hand, the $c_{al} = 1.5$
value allows no aligned edge.
This value is used as a reference, to show the performances of the other meshes compared ($c_{al} \leq 1$) with a pure uniform triangular mesh ($c_{al} = 1.5$).

Figure~\ref{fig:lshape:est:al:dof:15} illustrates the convergence of the relative a-posteriori error estimator
\begin{linenomath}
\begin{equation}
    \label{eq:eta:rel}
    \eta_{rel} := \frac{\eta_{\Omega}}{|||u|||_{\Omega}},
\end{equation} 
\end{linenomath}
and the energy norm relative error
\begin{linenomath}
\begin{equation}
    \label{eq:err:rel}
    |||e||| := \frac{|||U - \Pi^{\nabla}_k u|||_{\Omega}}{|||U|||_{\Omega}},
\end{equation}
\end{linenomath}
with respect to the number of DOFs, for different $c_{al}$ values, with VEM order $k = 3$.
Moreover, Figure~\ref{fig:lshape:est:al:slope:15} measures, for each iteration $m > 4$ of Process~\eqref{eq:schema}, the local convergence rate $\alpha$ based on the previous $(m-5, m-1)$ refinement iterations.
The analysis suggests that when $c_{\rho} > 1.0$ the parameter $c_{al}$ becomes less significant, as at convergence, all curves overlap and reach the optimal rate.
The results for orders $k=1$ and $k=2$ with different $c_{al}$ values are omitted as they exhibits a similar behaviour.
For completeness, Figure~\ref{fig:lshape:est:dof_est:15} displays the convergence of the error estimator for all the VEM orders with $c_{\rho} = 1.5$ and a fixed value of $c_{al} = 1.0$.

We conclude our analysis by exploring the scenario where $c_{\rho} \leq 1.0$.
In Figure~\ref{fig:lshape:est:al:dof:05} and Figure~\ref{fig:lshape:est:al:slope:05} we present the analysis of $c_{al}$ while fixing $k=3$ and $c_{\rho} = 0.5$.
We display again the convergence of relative estimator $\eta_{rel}$ of Equation~\eqref{eq:eta:rel} and the convergence rate $\alpha$ for each iteration $m$.
In Figure~\ref{fig:lshape:est:al:slope:05}, we report the rates of convergence and we observe that the best results are achieved with $c_{al} = 1.0$, specifically when $\# \mathcal{I}_E^e \leq 2$, $\forall e \in \mathcal{E}_E$, $\forall E \in \mathcal{T}^m_{\Omega}$.
Notably, the results obtained with $c_{al} = 1.0$ closely resemble those obtained with the pure triangular mesh of $c_{al} = 1.5$.
Conversely, when $c_{al} < 1.0$, the curves overlap and the convergence of the estimator appears sub-optimal.
We believe that this sub-optimal behaviour may be attributed to the presence of excessively small edges in the mesh.
In Figure~\ref{fig:lshape:est:al:qE:05} and Figure~\ref{fig:lshape:est:al:qE:10} we provide a statistical analysis of the quantities $h_E$ and $r_E$, for $c_{al} = 0.5$ and $c_{al} = 1.0$, respectively.
The quantity 
\begin{linenomath}
\begin{equation}
    \label{eq:rhoE}
    \rho_E := \min \{h_E, r_E\},\ \forall E \in \mathcal{T}_{\Omega}^m.
\end{equation}
\end{linenomath}
is the one used in the \texttt{CHECK-QUALITY} Algorithm~\ref{alg:quality}.
We present the first and third quartile values (q1-q3 area) and the average values (avg curve).
The analysis of the statistics reveals that when the convergence rates $\alpha$ are not optimal ($c_{al} = 0.5$), the average curve of $\rho_E$ deviates from the average curve $r_E$, and the value of $h_E$ is nearly equal or lower than $r_E$.
This discrepancy indicates that, on average, when the length of the minimum edge of the mesh elements ($h_E$) is not significantly higher than the inner radius ($r_E$), the VEM assumptions regarding the quality of the element edges are not fulfilled, jeopardizing the optimal convergence of the numeric solution, as discussed in Section~\ref{sec:setting}.
For these reasons, we select $c_{al} = 1.0$ as the optimal parameter value, and we present the convergence of the estimator $\eta_{rel}$ in Figure~\ref{fig:lshape:est:dof_est:05} for all VEM orders $k=1$, $k=2$, and $k=3$.

Before proceeding with the DFN tests, we stress that we obtain the optimal order in the VEM approximations selected and for both choices of $c_{\rho}$, thanks to the introduction of the parameter $c_{al}$ in Check~\eqref{check:two}, despite the differences in the resulting mesh shapes.

\subsection{Test 2: Regular DFN with three Fractures}
\begin{figure}[!h]
    \centering
    \begin{subfigure}[]{0.32\textwidth}
        \centering
        \includegraphics[width=\textwidth]{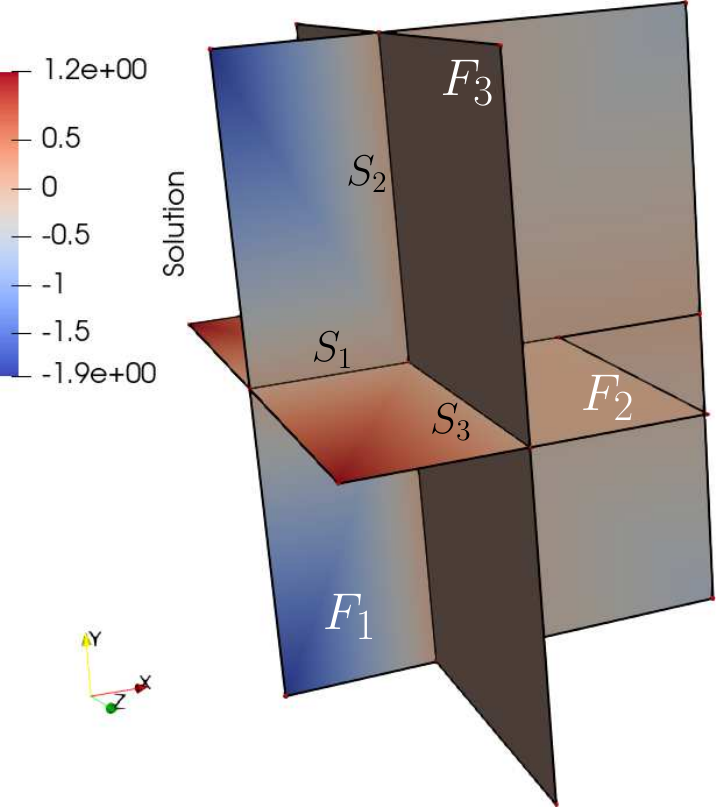}
        \caption{DFN $\mathcal{T}^0_{\Omega}$}
        \label{fig:tipetut:network:t0:05}
    \end{subfigure}
    \begin{subfigure}[]{0.32\textwidth}
        \centering
        \includegraphics[width=\textwidth]{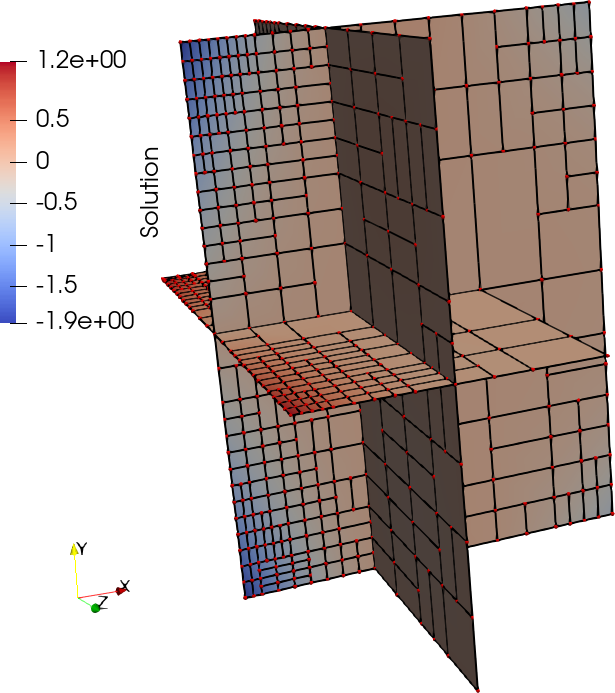}
        \caption{Step $m=15$, $c_{\rho} = 0.5$}
        \label{fig:tipetut:network:t15:05}
    \end{subfigure}
    \begin{subfigure}[]{0.32\textwidth}
        \centering
        \includegraphics[width=\textwidth]{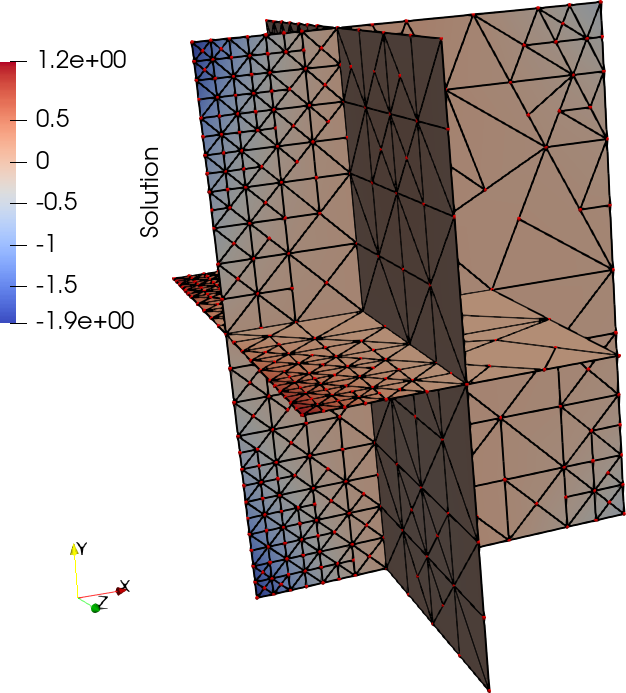}
        \caption{Step $m=15$, $c_{\rho} = 1.5$}
        \label{fig:tipetut:network:t15:15}
    \end{subfigure}
    \caption{Test $2$ - Network with solution, $k = 1$, $c_{al} = 1.0$.}
    \label{fig:tipetut:network}
\end{figure}
\begin{figure}[!h]
    \centering
    \begin{subfigure}[]{0.32\textwidth}
        \centering
        \includegraphics[width=\textwidth]{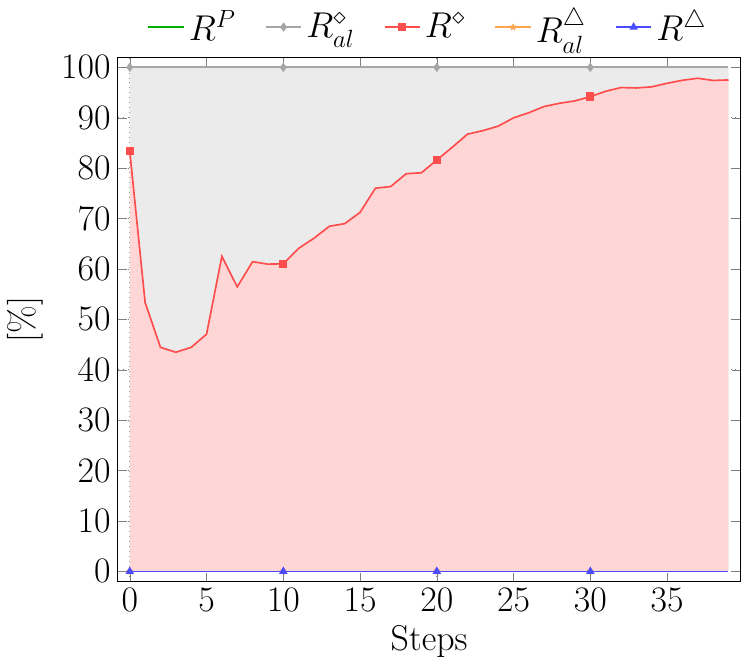}
        \caption{$c_\rho = 0.5$}
        \label{fig:tipetut:m_quality:tri:05}
    \end{subfigure}
    \begin{subfigure}[]{0.32\textwidth}
        \centering
        \includegraphics[width=\textwidth]{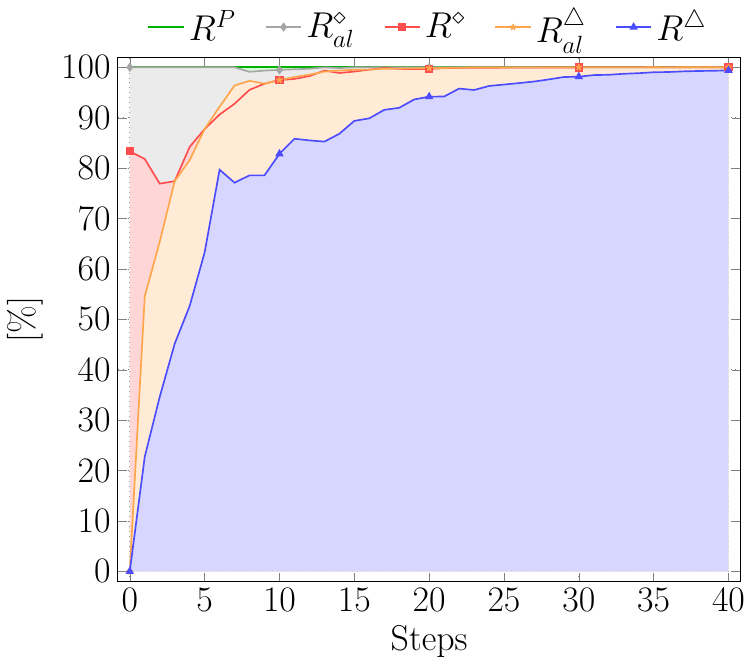}
        \caption{$c_\rho = 1.5$}
        \label{fig:tipetut:m_quality:tri:15}
    \end{subfigure}
    \caption{Test $2$ - $R^{\star}$ versus $m$-steps, $k = 1$, $c_{al} = 1.0$.}
    \label{fig:tipetut:dofs_ar}
\end{figure}
\begin{figure}[!h]
    \centering
    \begin{subfigure}[]{0.4\textwidth}
        \centering
        \includegraphics[width=\textwidth]{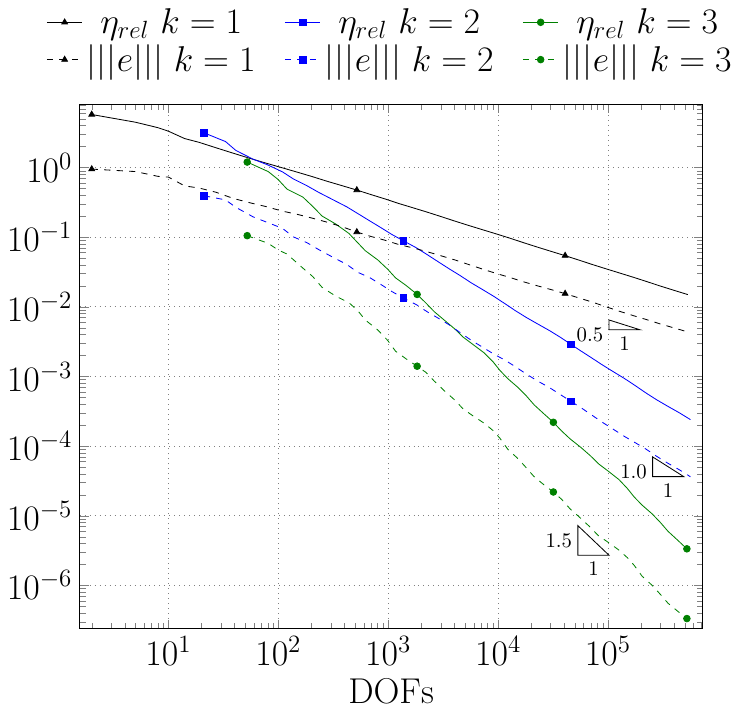}
        \caption{$c_\rho = 0.5$}
        \label{fig:tipetut:est:dof_est:05}
    \end{subfigure}
    \begin{subfigure}[]{0.4\textwidth}
        \centering
        \includegraphics[width=\textwidth]{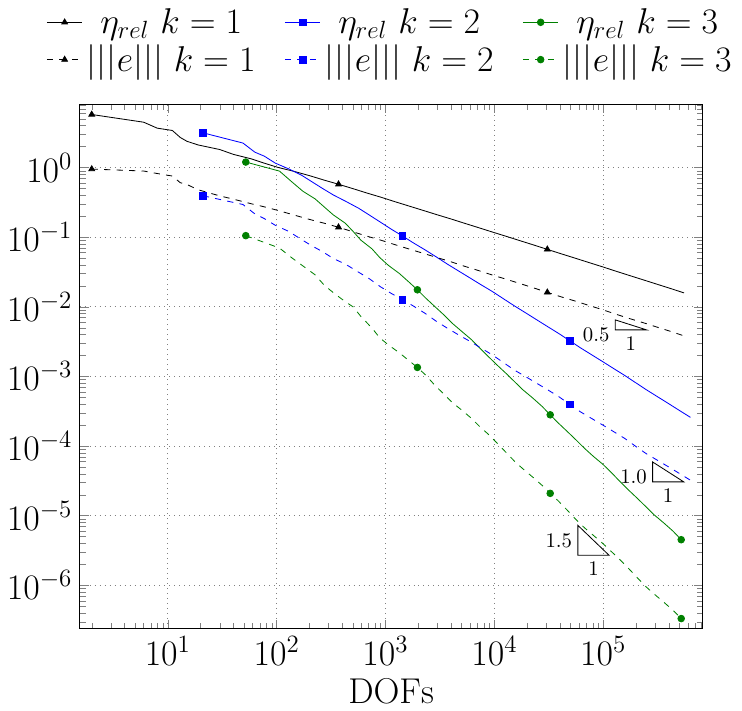}
        \caption{$c_\rho = 1.5$}
        \label{fig:tipetut:est:dof_est:15}
    \end{subfigure}
    \caption{Test $2$ - Relative Estimator ($\eta_{rel}$) and energy error ($|||e|||$) versus DOFs, $c_{al} = 1.0$.}
    \label{fig:tipetut:est}
\end{figure}
\begin{table}[!h]
    \caption{Test $2$ - Rates of convergence $\alpha$ for the estimator $\eta_{rel}$ and the error $e_{H^1}$, $c_{al} = 1.0$.}
    \label{tab:tipetut:est}
    \centering
    \includegraphics[width=0.45\textwidth]{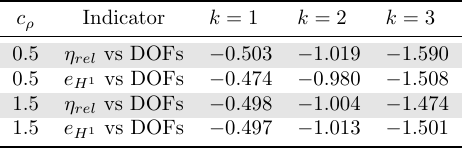}
\end{table}
We present a benchmark for a three-fractures DFN test.
The detailed problem information are available in \cite{BENEDETTO2016148}.
We prescribe the exact solution on each fracture as follows:
\begin{itemize}
    \item[$F_1$] $-\frac{1}{10} (x + \frac{1}{2})\left[8 x y (x^2+y^2)\arctan{(y, x)}+x^3\right]$,
    \item[$F_2$] $-\frac{1}{10} (x + \frac{1}{2})x^3(1-8\pi |z|)$,
    \item[$F_3$] $y(y - 1)(y + 1)(z - 1) z$.
\end{itemize}
Figure~\ref{fig:tipetut:network:t0:05} illustrates the network with the initial mesh $\mathcal{T}_{\Omega}^0$.
As in the previous case, the subsequent analysis employs VEM order $k= \{1, 2, 3\}$, $c_{\rho} \in \{0.5, 1.5\}$, and $c_{al} = 1.0$.
We do not report any analysis varying $c_{al}$ due to negligible differences obtained, attributed to the large regularity of the original network geometry.
This test, in contrast to the concave L-shape domain, encompasses all convex domains with orthogonal intersections.
The iterative scheme~\eqref{eq:schema} terminates when the total number of DOFs reaches $5 \cdot 10^5$.

Figure~\ref{fig:tipetut:network:t15:05} and Figure~\ref{fig:tipetut:network:t15:15} depict the mesh and solution at refinement step $m=15$ for $c_{\rho} = 0.5$ and $c_{\rho} = 1.5$, respectively.
Significant differences in the generated meshes are evident.
Indeed, the analysis of the $R^{\star}$ values in Figure~\ref{fig:tipetut:m_quality:tri:05} reveals that for $c_{\rho} < 1$ the mesh, $\mathcal{T}^m_{\Omega}$, converges to a full quadrilateral mesh.
This behaviour is attributed, once again, to the regularity of the original network geometry.
On the other hand, in Figure~\ref{fig:tipetut:m_quality:tri:15}, when $c_{\rho} > 1$, the mesh $\mathcal{T}^m_{\Omega}$ approaches to a fully triangular discretization.
The plots of the $AR^{\star}_E$ indicators are omitted since all indicators trivially converge to the characteristic values of the reference square $[0,1]\times[0,1]$ and the reference triangle shown in Figure~\ref{fig:lshape:reference} for $c_{\rho} = 0.5$ and $c_{\rho} = 1.5$, respectively.

It is important to note that the shape regularity observed in this network is atypical in real DFN applications.
Nevertheless, this test was specifically designed to have a known exact solution.

Figures~\ref{fig:tipetut:est} present the convergence of the relative estimator $\eta_{rel}$ of Equation~\eqref{eq:eta:rel} and the energy norm relative error of Equation~\eqref{eq:err:rel} for each DOFs.
Additionally, Table~\ref{tab:tipetut:est} details the convergence rates $\alpha$, computed on the last $5$ iterations of refinement.
It can be affirmed that optimal convergence rates $\alpha$ are achieved for both the estimator $\eta_{rel}$ and the relative error $|||e|||$ for each VEM order $k$ and each choice of $c_{\rho}$.

\subsection{Test 3: Random DFN}
\begin{figure}[!h]
    \centering
    \begin{subfigure}[]{0.32\textwidth}
        \centering
        \includegraphics[width=\textwidth]{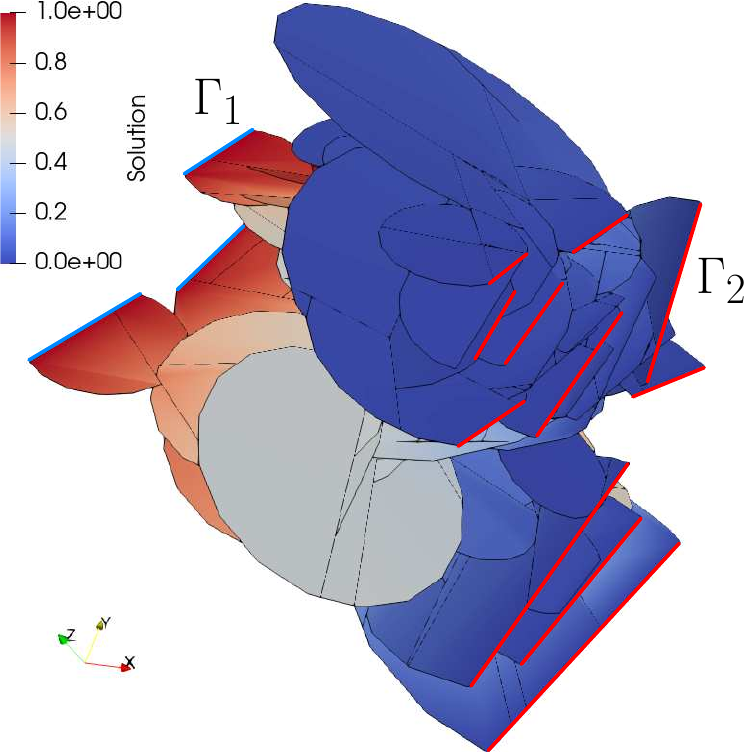}
        \caption{$\Omega$ - $\mathcal{T}^0_{\Omega}$}
        \label{fig:frac86:network:t0:05}
    \end{subfigure}
    \begin{subfigure}[]{0.32\textwidth}
        \centering
        \includegraphics[width=\textwidth]{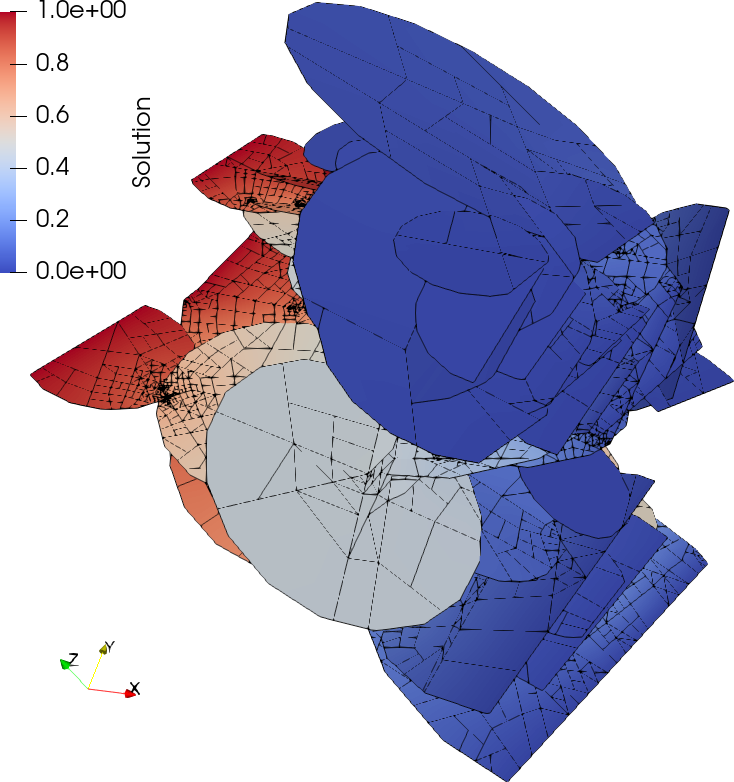}
        \caption{$\Omega$ - $m=15$, $c_\rho = 0.5$}
        \label{fig:frac86:network:t15:05}
    \end{subfigure}
    \begin{subfigure}[]{0.32\textwidth}
        \centering
        \includegraphics[width=\textwidth]{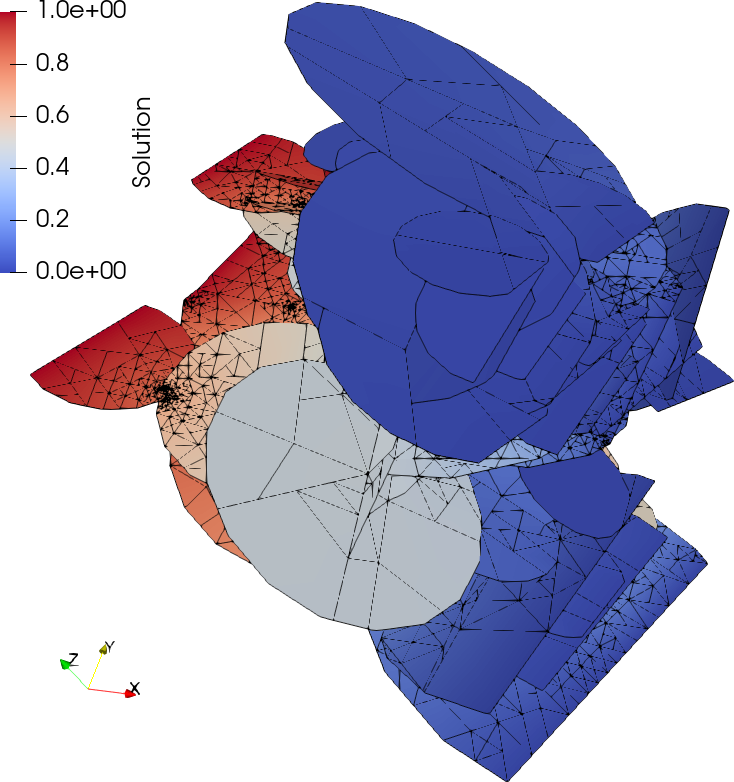}
        \caption{$\Omega$ - $m=15$, $c_\rho = 1.5$}
        \label{fig:frac86:network:t15:15}
    \end{subfigure}\\
    \begin{subfigure}[]{0.32\textwidth}
        \centering
        \includegraphics[width=\textwidth]{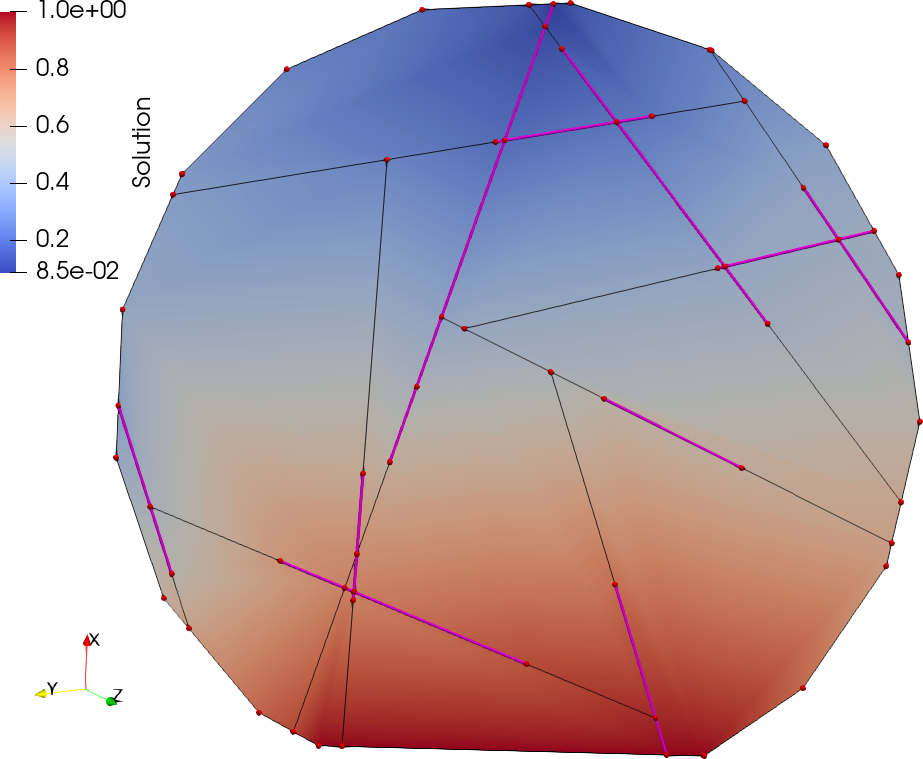}
        \caption{$F_{72}$ - $\mathcal{T}^0_{\Omega}$}
        \label{fig:frac86:F71:t0:05}
    \end{subfigure}
    \begin{subfigure}[]{0.32\textwidth}
        \centering
        \includegraphics[width=\textwidth]{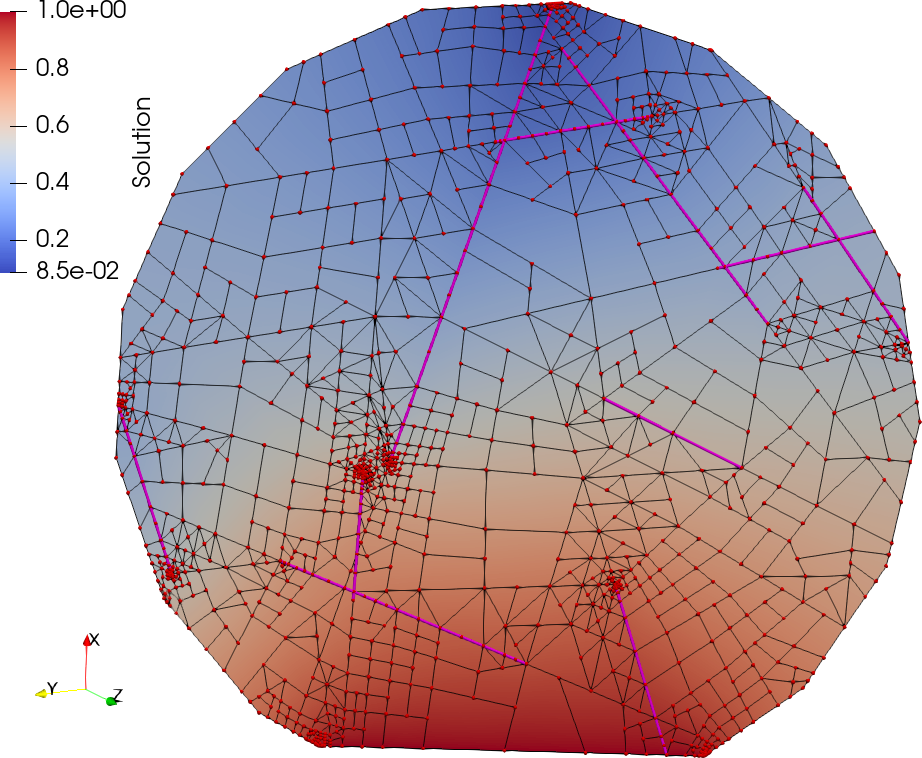}
        \caption{$F_{72}$ - $m=15$, $c_\rho = 0.5$}
        \label{fig:frac86:F71:t15:05}
    \end{subfigure}
    \begin{subfigure}[]{0.32\textwidth}
        \centering
        \includegraphics[width=\textwidth]{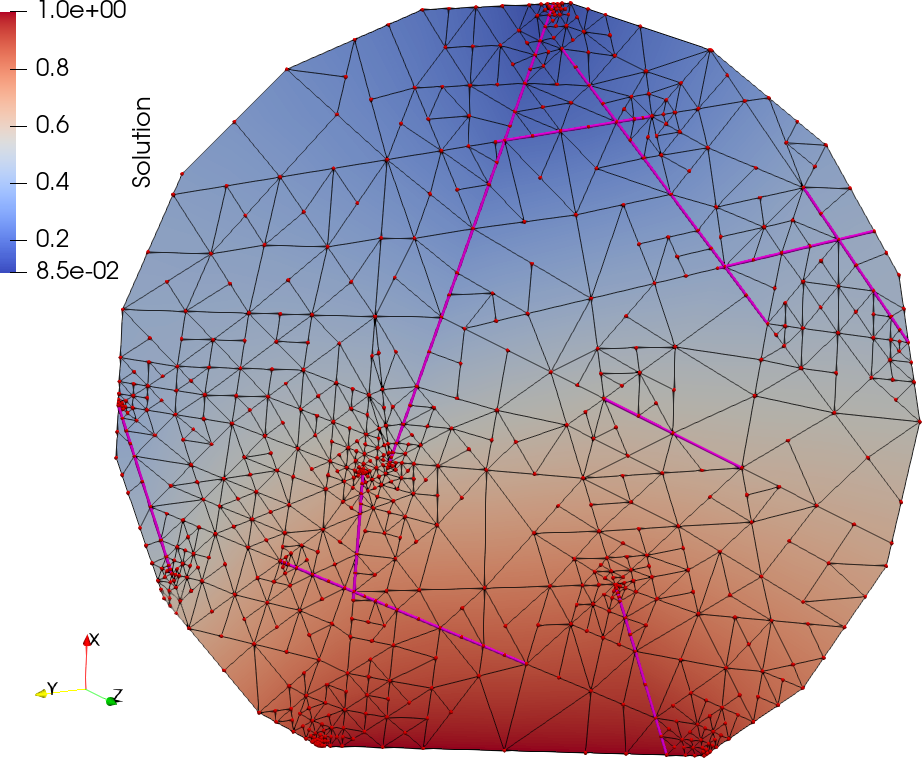}
        \caption{$F_{72}$ - $m=15$, $c_\rho = 1.5$}
        \label{fig:frac86:F71:t15:15}
    \end{subfigure}
    \caption{Test $3$ - Solution on the network $\Omega$ and on $F_{72}$, $k=1$, $c_{al} = 1.0$.}
    \label{fig:frac86:network}
\end{figure}
\begin{figure}[!h]
    \centering
    \begin{subfigure}[]{0.32\textwidth}
        \centering
        \includegraphics[width=\textwidth]{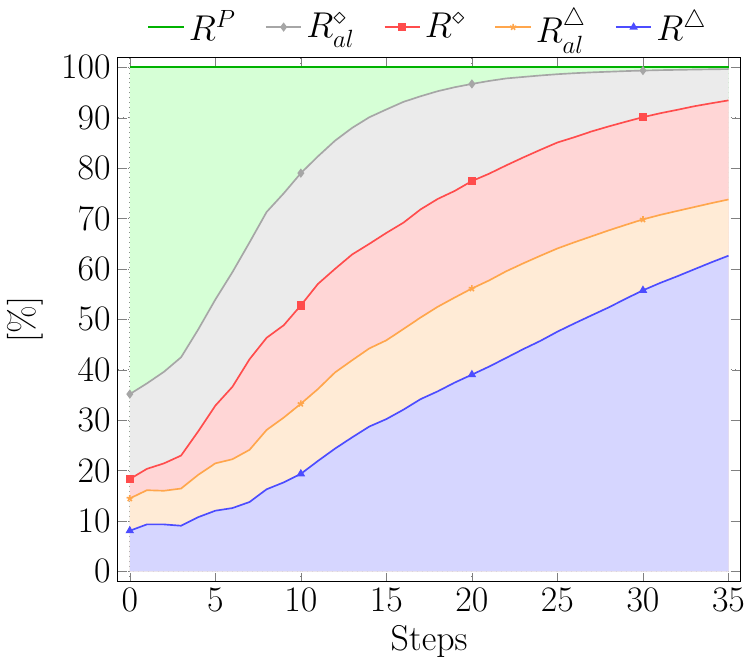}
        \caption{$R^{\star}$ vs $m$-steps, $c_\rho = 0.5$}
        \label{fig:frac86:m_quality:rho:05}
    \end{subfigure}
    \begin{subfigure}[]{0.32\textwidth}
        \centering
        \includegraphics[width=\textwidth]{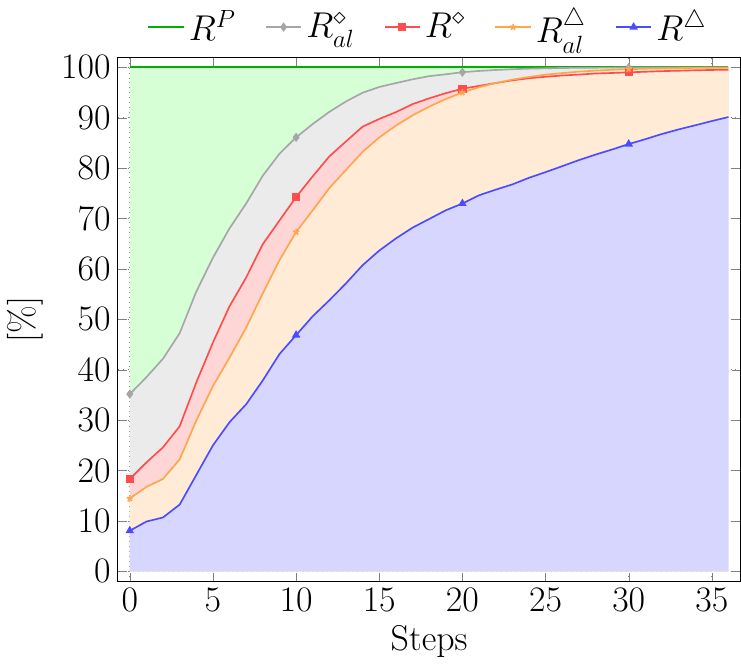}
        \caption{$R^{\star}$ vs $m$-steps, $c_\rho = 1.5$}
        \label{fig:frac86:m_quality:rho:15}
    \end{subfigure}
    \begin{subfigure}[]{0.32\textwidth}
        \centering
        \includegraphics[width=\textwidth]{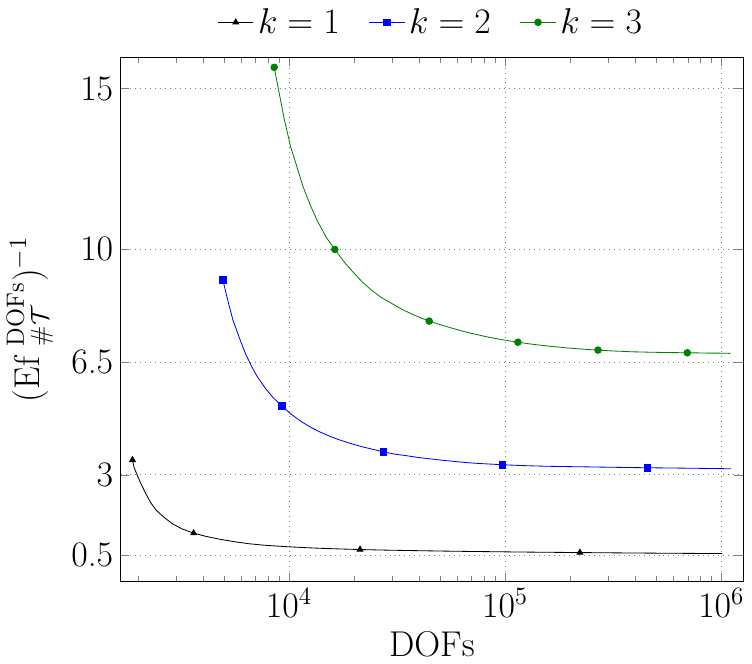}
        \caption{$(\text{Ef }^{\text{DOFs}}_{\# \mathcal{T}})^{-1}$, $c_\rho = 1.5$ vs DOFs}
        \label{fig:frac86:m_quality:ePT:15}
    \end{subfigure}
    \caption{Test $3$ - Mesh $\mathcal{T}^m_{\Omega}$ shapes, $k=1$, $c_{al} = 1.0$.}
    \label{fig:frac86:m_quality}
\end{figure}
\begin{figure}[!h]
    \centering
    \begin{subfigure}[]{0.32\textwidth}
        \centering
        \includegraphics[width=\textwidth]{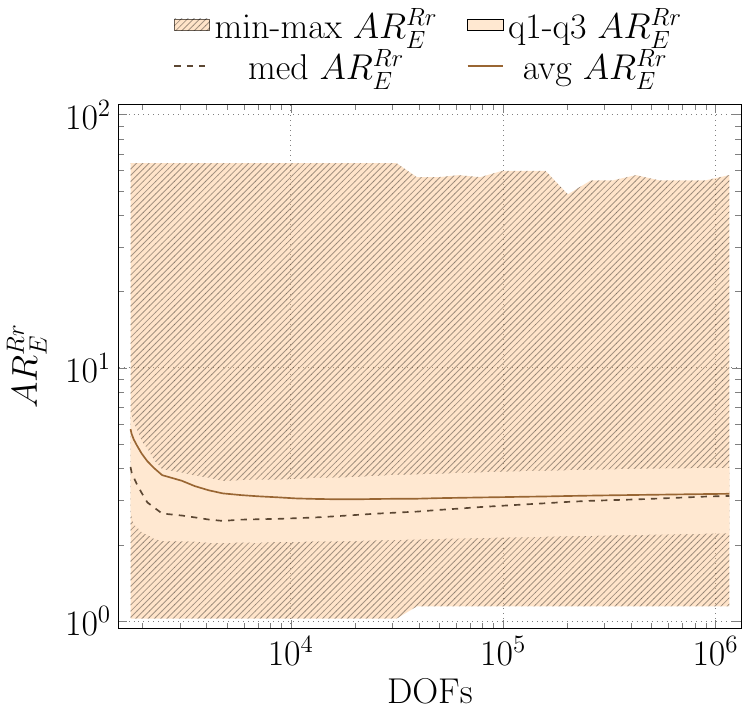}
        \caption{$\Omega$ - $AR^{Rr}_E$ vs DOFs, $c_{\rho} = 0.5$}
        \label{fig:frac86:dofs_ar:Rr:05}
    \end{subfigure}
    \begin{subfigure}[]{0.32\textwidth}
        \centering
        \includegraphics[width=\textwidth]{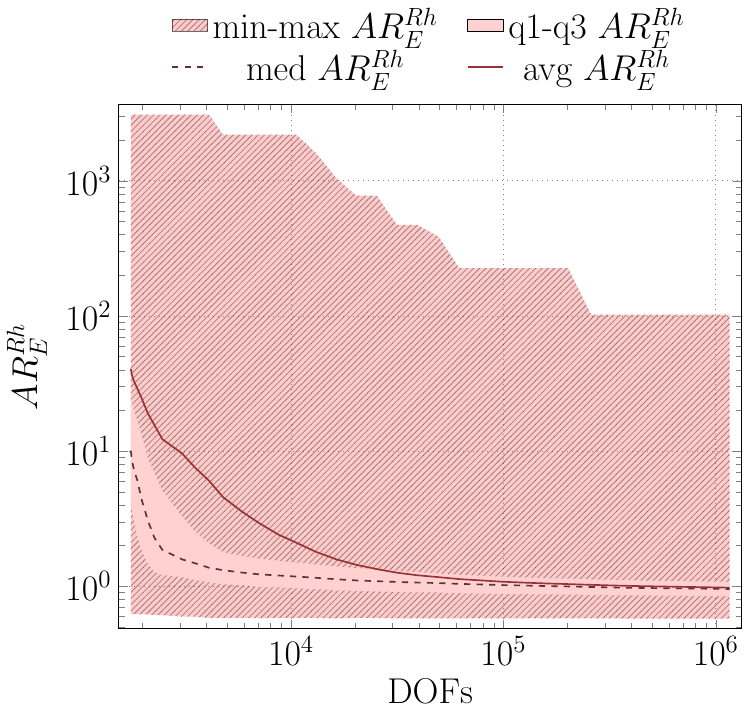}
        \caption{$\Omega$ - $AR^{Rh}_E$ vs DOFs, $c_{\rho} = 0.5$}
        \label{fig:frac86:dofs_ar:Rh:05}
    \end{subfigure}\\
    \begin{subfigure}[]{0.32\textwidth}
        \centering
        \includegraphics[width=\textwidth]{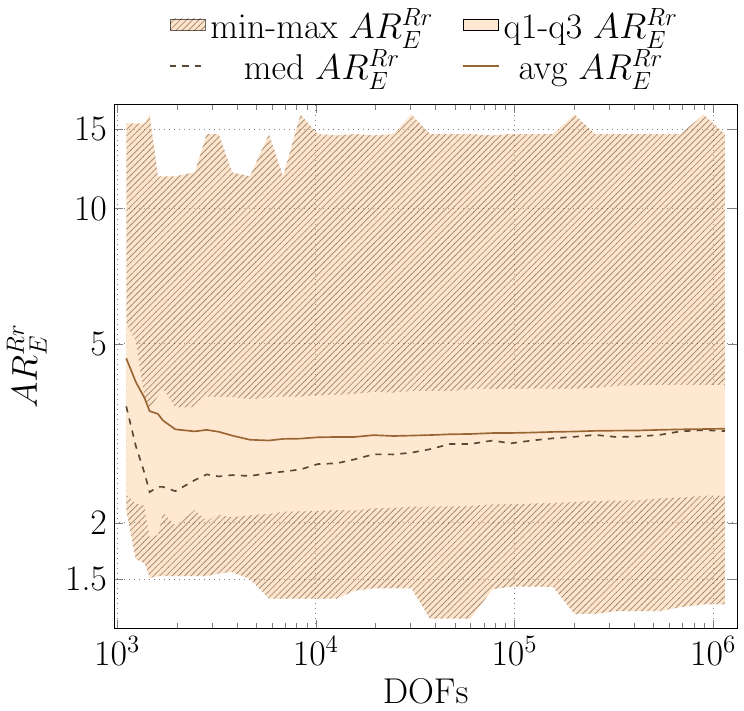}
        \caption{$F_{72}$ - $AR^{Rr}_E$ vs DOFs, $c_{\rho} = 0.5$}
        \label{fig:frac86:F71:dofs_ar:Rr:05}
    \end{subfigure}
    \begin{subfigure}[]{0.32\textwidth}
        \centering
        \includegraphics[width=\textwidth]{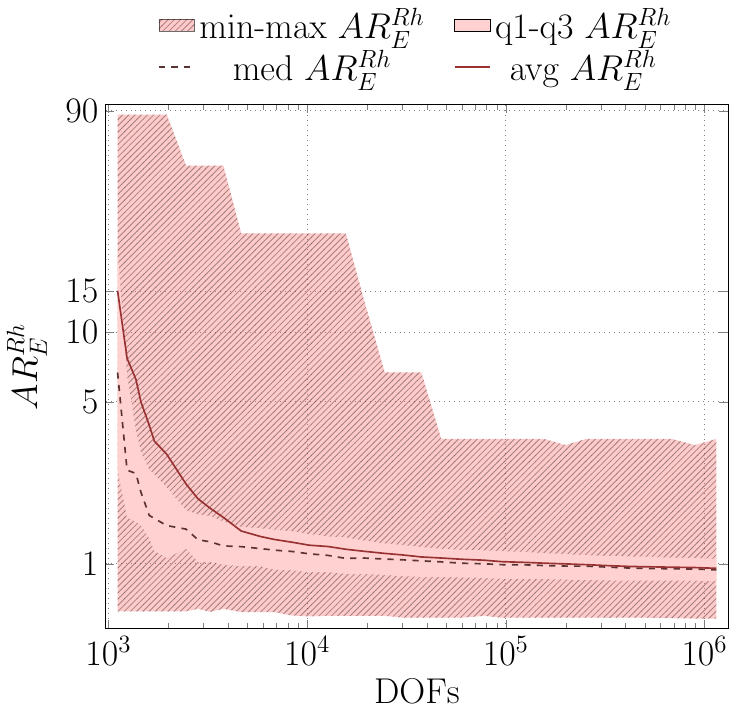}
        \caption{$F_{72}$ - $AR^{Rh}_E$ vs DOFs, $c_{\rho} = 0.5$}
        \label{fig:frac86:F71:dofs_ar:Rh:05}
    \end{subfigure}
    \caption{Test $3$ - Mesh $\mathcal{T}^m_{\Omega}$ quality, $k=1$, $c_{al} = 1.0$.}
    \label{fig:frac86:dofs_ar}
\end{figure}
\begin{figure}[!h]
    \centering
    \begin{subfigure}[]{0.38\textwidth}
        \centering
        \includegraphics[width=\textwidth]{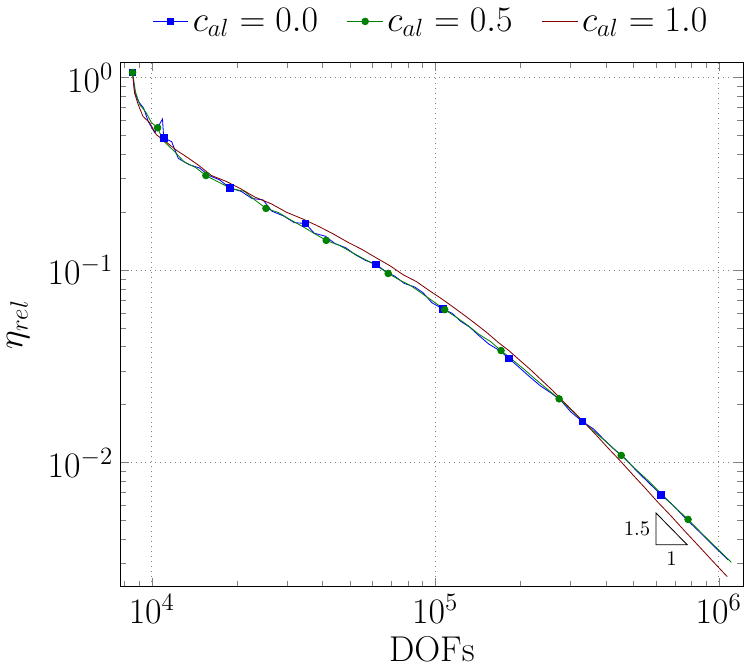}
        \caption{Estimator vs DOFs}
        \label{fig:frac86:est:al:dof:05}
    \end{subfigure}
    \begin{subfigure}[]{0.38\textwidth}
        \centering
        \includegraphics[width=\textwidth]{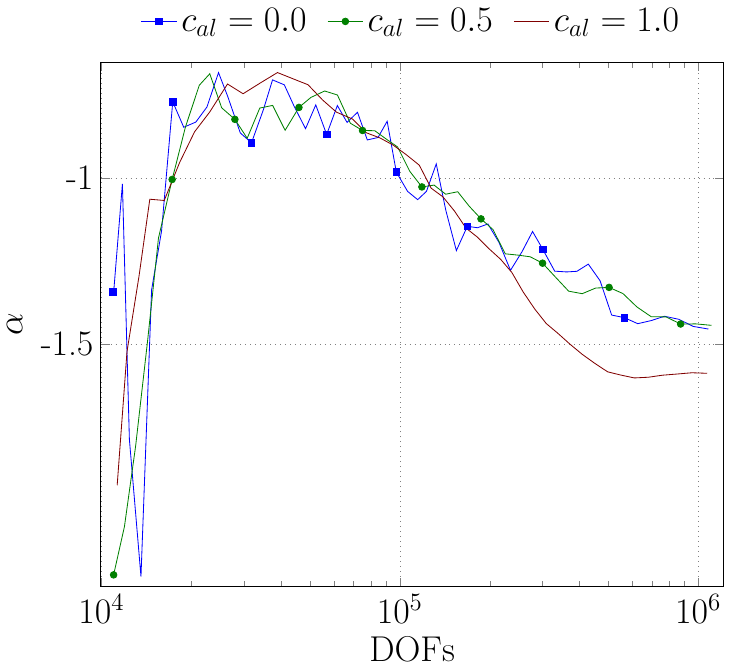}
        \caption{Slope $\alpha$ vs DOFs}
        \label{fig:frac86:est:al:slope:05}
    \end{subfigure}
    \begin{subfigure}[]{0.38\textwidth}
        \centering
        \includegraphics[width=\textwidth]{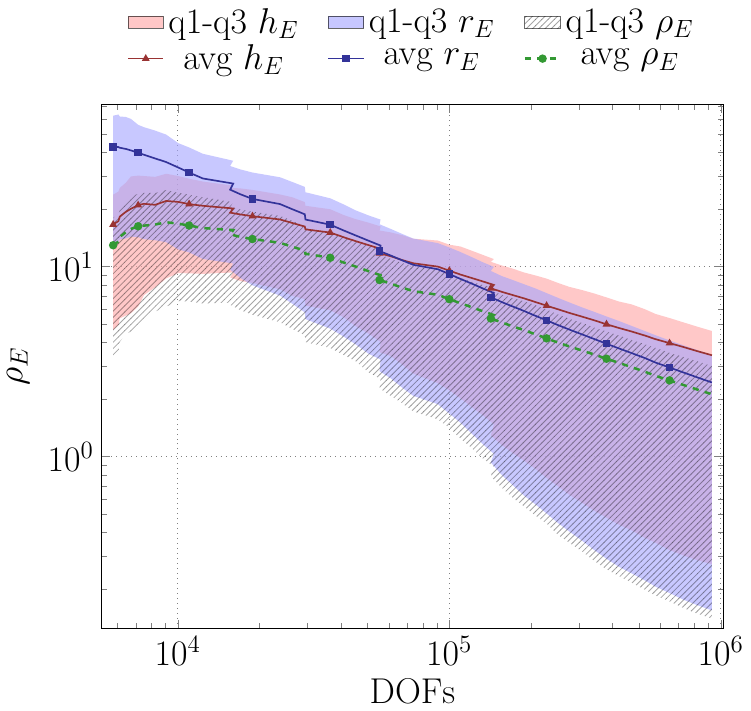}
        \caption{$\rho_E$ vs DOFs, $c_{al} = 0.5$}
        \label{fig:frac86:est:al:q:05}
    \end{subfigure}
    \begin{subfigure}[]{0.38\textwidth}
        \centering
        \includegraphics[width=\textwidth]{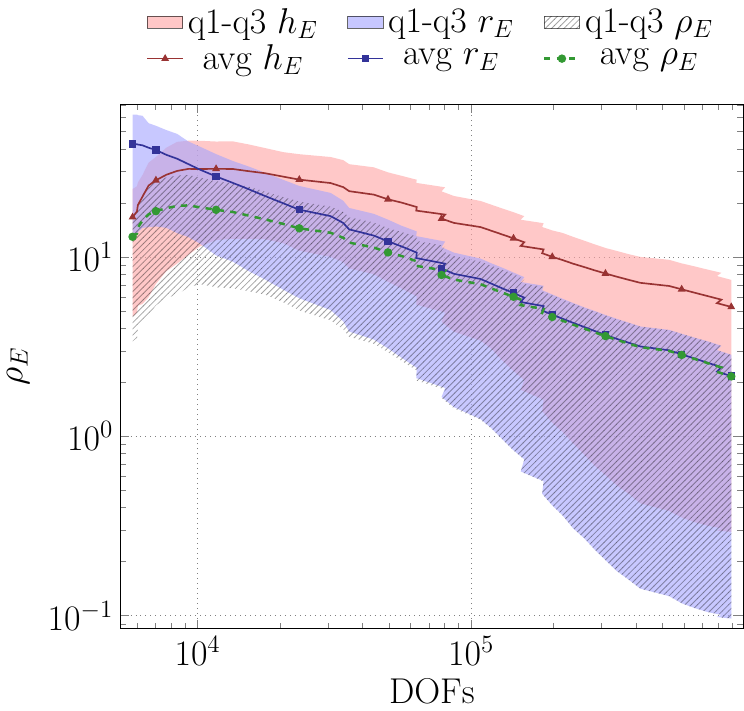}
        \caption{$\rho_E$ vs DOFs, $c_{al} = 1.0$}
        \label{fig:frac86:est:al:q:10}
    \end{subfigure}
    \caption{Test $3$ - $c_{al}$ analysis, $c_{\rho} = 0.5$, $k=3$.}
    \label{fig:frac86:est:al}
\end{figure}
\begin{figure}[!h]
    \centering
    \begin{subfigure}[]{0.38\textwidth}
        \centering
        \includegraphics[width=\textwidth]{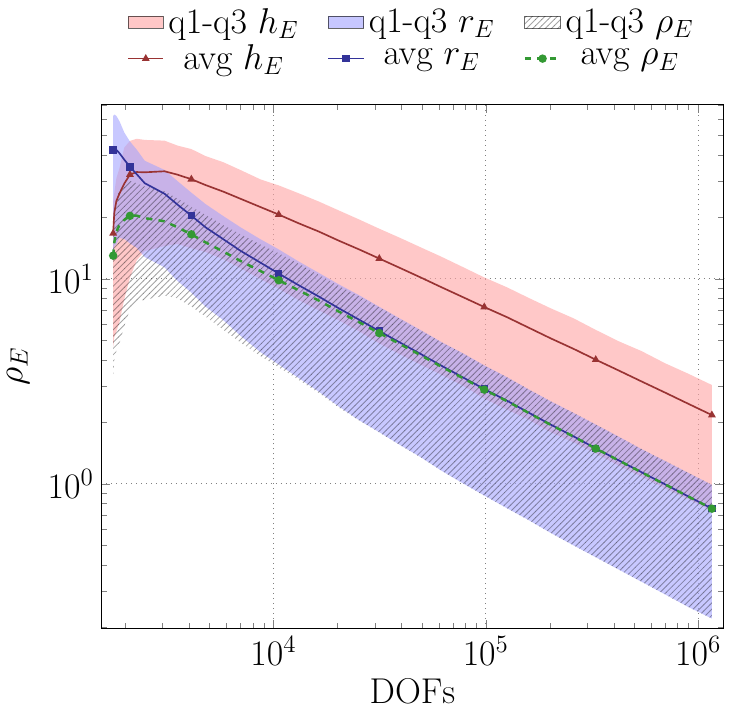}
        \caption{$\rho_E$ vs DOFs, $c_\rho = 0.5$, $k=1$}
        \label{fig:frac86:m_rhoE:O1:05}
    \end{subfigure}
    \begin{subfigure}[]{0.38\textwidth}
        \centering
        \includegraphics[width=\textwidth]{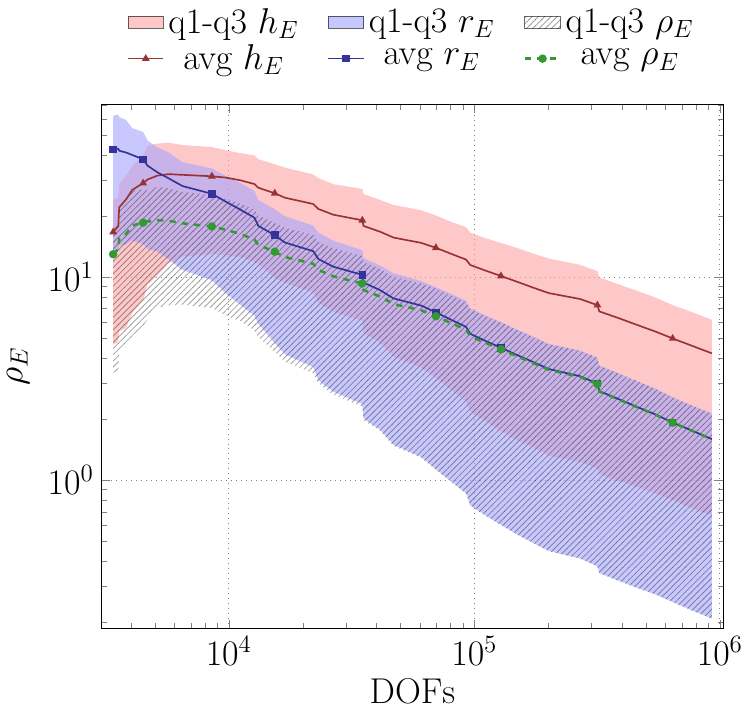}
        \caption{$\rho_E$ vs DOFs, $c_\rho = 0.5$, $k=2$}
        \label{fig:frac86:m_rhoE:O2:05}
    \end{subfigure}
    \begin{subfigure}[]{0.38\textwidth}
        \centering
        \includegraphics[width=\textwidth]{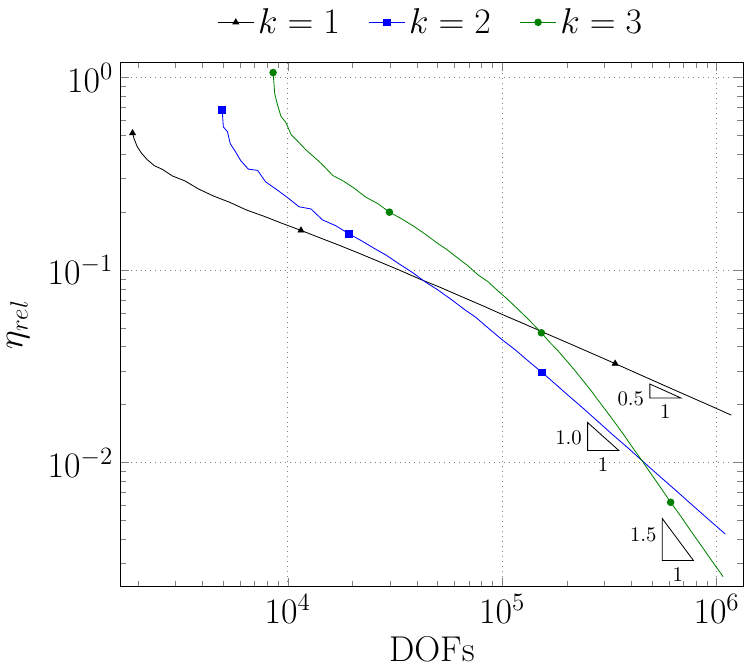}
        \caption{$\eta_{rel}$ vs DOFs, $c_\rho = 0.5$}
        \label{fig:frac86:est:dof_est:05}
    \end{subfigure}
    \begin{subfigure}[]{0.38\textwidth}
        \centering
        \includegraphics[width=\textwidth]{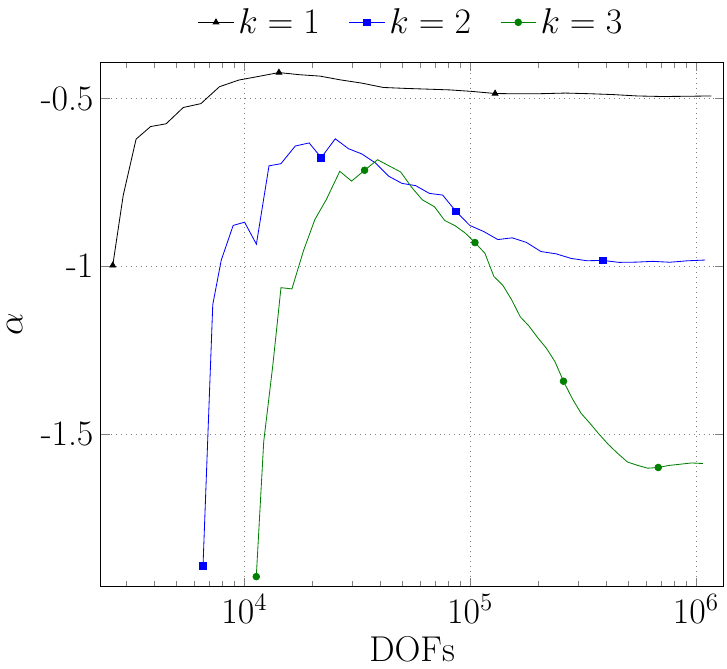}
        \caption{$\alpha$ vs DOFs, $c_\rho = 0.5$}
        \label{fig:frac86:est:dof_slope:05}
    \end{subfigure}
    \caption{Test $3$ - Estimator $\eta_{rel}$ analysis, $c_{al} = 1.0$.}
    \label{fig:frac86:est}
\end{figure}
\begin{figure}[!h]
    \centering
    \begin{subfigure}[]{0.4\textwidth}
        \centering
        \includegraphics[width=\textwidth]{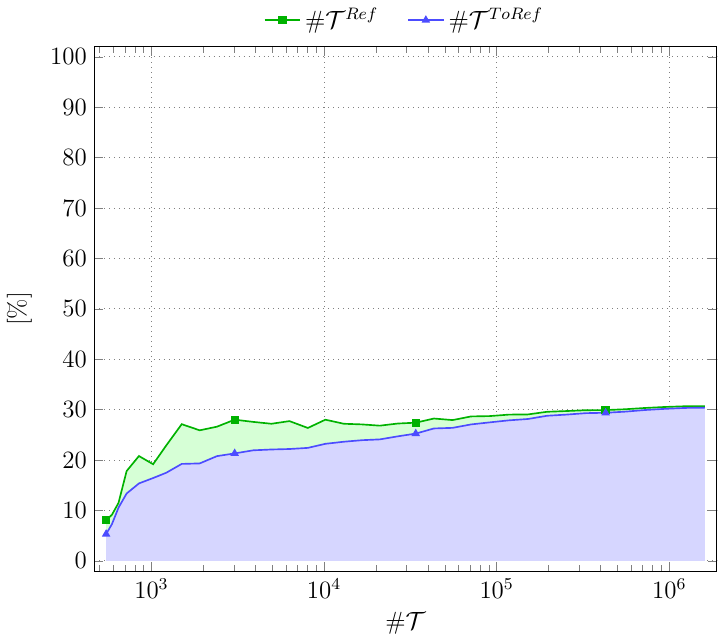}
        \caption{$\# \mathcal{T}^{\star}_{\Omega}$ vs $\# \mathcal{T}_{\Omega}$}
        \label{fig:frac86:ref:cells:05}
    \end{subfigure}
    \begin{subfigure}[]{0.4\textwidth}
        \centering
        \includegraphics[width=\textwidth]{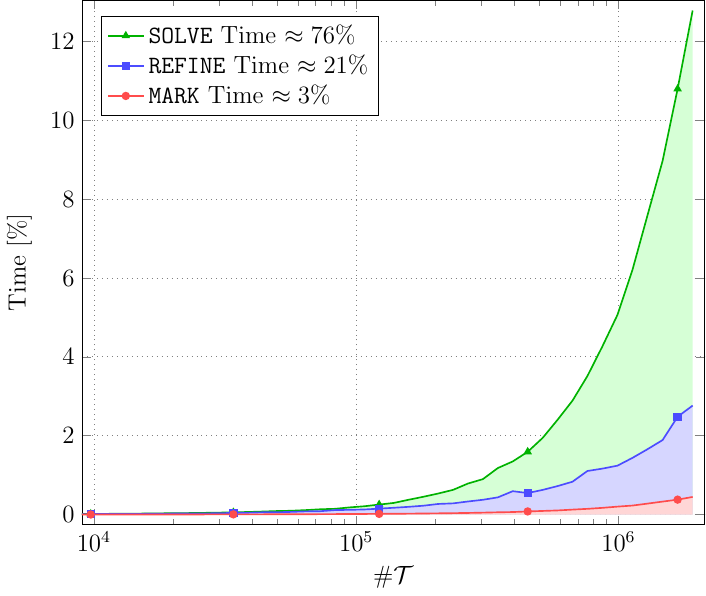}
        \caption{Time $[\%]$ vs $\# \mathcal{T}_{\Omega}$}
        \label{fig:frac86:ref:time:05}
    \end{subfigure}
    \caption{Test $3$ - Refinement Impact analysis, $c_\rho = 0.5$, $c_{al} = 1.0$, $k=1$.}
    \label{fig:frac86:ref}
\end{figure}
To conclude the analysis, we discuss the results obtained for a more realistic DFN case.
We select the benchmark test analysed in \cite{BERRONE2021103502}.
The network, illustrated in Figure~\ref{fig:frac86:network:t0:05}, features a computed \emph{minimal mesh} $\mathcal{T}^0_{\Omega}$ generated by $86$ fractures and $159$ fracture intersections.
We compute the solution to Problem~\eqref{eq:prob:cont} by imposing two constant Dirichlet conditions, $\Gamma_1 = 1.0$ and $\Gamma_2 = 0.0$, on the fracture edges that intersect the planes $x=0.0$ and $x=1000.0$, respectively.
Other fracture edges have homogeneous Neumann conditions.
With blue and red thick edges in Figure~\ref{fig:frac86:network:t0:05} we highlight $\Gamma_1$ and $\Gamma_2$, respectively.

This DFN test is generated with random fracture orientation, size, and localization drawn from random distributions computed from a real fractured media.
Consequently, the initial network mesh $\mathcal{T}^0_{\Omega}$ consists of several generic polygonal elements with differing size and a notable number of aligned edges.
Figure~\ref{fig:frac86:F71:t0:05} offers a closer look at the initial mesh of fracture $F_{72} \in \Omega$ as a detail of the global mesh $\mathcal{T}^0_{\Omega}$.
The fracture intersections are highlighted with purple lines.
The complexity and the poor element quality are recognisable.

The analysis of the mesh involves VEM order $k \in \{1, 2, 3\}$, $c_{\rho} \in \{0.5, 1.5\}$, and $c_{al} = \{0.0, 0.5, 1.0\}$.
We stop the iterative scheme~\eqref{eq:schema} when the total number of DOFs reaches $10^6$.

Figures~\ref{fig:frac86:network:t15:05}-\ref{fig:frac86:network:t15:15} and Figures~\ref{fig:frac86:F71:t15:05}-\ref{fig:frac86:F71:t15:15} illustrate the mesh obtained after $15$ steps in the network $\Omega$ and in $F_{72}$ for $c_{\rho} = 0.5$ and $c_{\rho} = 1.5$, respectively.
Notably, the mesh refinement algorithm strategically focuses around the tips of fracture intersections where the solution exhibits a singularity $H^{\frac{3}{2}-\epsilon}$.

In the plots of Figure~\ref{fig:frac86:m_quality:rho:05} and Figure~\ref{fig:frac86:m_quality:rho:15}, we measure the $R^{\star}$ indicators for both the $c_{\rho}$ values.
The fracture mesh detail and the analysis of the plots revel that higher values of $c_{\rho}$ lead to an increased number of mesh triangles,  consistent with the observations from previous tests.
When $c_{\rho} = 1.5$, almost the entire number of elements consist of triangular or nearly triangular ($R^{\Delta} + R^{\Delta}_{al} \approx 100 \%$).
In addition, despite the initially prevalence of generic polygons in $\mathcal{T}^0_{\Omega}$ ($R^P \approx 65 \%$), both meshes ultimately converge to tessellations containing predominantly triangular or quadrilateral elements or those that are nearly triangular and quadrilateral ($R^P \approx 0$).
This highlights the robustness of the approach in improving the shape-uniformity of the elements on generic polygonal meshes.

Figure~\ref{fig:frac86:m_quality:ePT:15} presents the inverse of $\text{Ef }^{\text{DOFs}}_{\# \mathcal{T}_{\Omega}}$ indicator against the number of DOFs for $c_{\rho} = 1.5$.
The data for $c_{\rho} = 0.5$ is omitted as the plots are qualitatively similar.
This indicator $(\text{Ef }^{\text{DOFs}}_{\# \mathcal{T}_{\Omega}})^{-1}$ is presented as it converges to known values for triangular meshes and it evaluates the inefficiency of the mesh elements in approximating the numerical solution as the refinement process progresses.
Indeed, the indicator starts from very high values and rapidly converges to characteristic values of a triangular mesh, namely $0.5$, $3$ and $6.5$ for $k=1$, $k=2$ and $k=3$, respectively.
Thus, the measurement of $(\text{Ef }^{\text{DOFs}}_{\# \mathcal{T}_{\Omega}})^{-1}$ ensures that, by the end of the refinement process, maximum mesh efficiency in terms of DOFs is achieved.

In Figures~\ref{fig:frac86:dofs_ar}, we present a comprehensive statistical analysis of the $AR^{\star}_E$ quantities, where $\star \in \{Rh, Rr\}$, for the $c_{\rho} = 0.5$ case in the DFN $\Omega$ and in the fracture $F_{72}$.
The test with $c_{\rho} = 1.5$ yields similar results, so we have chosen to omit the corresponding plots.
For each step $m$, we measure the minimum and maximum values (min-max $AR_E^{\star}$ area), the first and third quartile values (q1-q3 $AR_E^{\star}$ area), the median value (med $AR_E^{\star}$ curve), and the average value (avg $AR_E^{\star}$ curve).
In the DFN, Figures~\ref{fig:frac86:dofs_ar:Rr:05}-\ref{fig:frac86:dofs_ar:Rh:05} show that as the number of DOFs increases, all the average curves approach the median curves.
This indicates an improvement in the quality of mesh elements at convergence.
The maximum values of the DFN $AR^{Rr}_E$ remain stable during the iterations because some fractures are almost untouched by the refinement process, where the solution is nearly constant (dead-end fractures).
It's important to note that these poor-quality elements are in a limited number because the median and the mean curves are order of magnitude far from the maximum values of the plot.
On the other hand, in the region where the refinement process is exploited, the effect of the mesh quality improvement are visible.
Indeed, for the local tessellation of fracture $F_{72}$, Figures~\ref{fig:frac86:F71:dofs_ar:Rr:05}-\ref{fig:frac86:F71:dofs_ar:Rh:05} reveal that particularly noteworthy is the decay of the $AR_E^{Rh}$ indicator, starting from a value of $90$ and ending at a small maximum value lower than $5$.

In Figures~\ref{fig:frac86:est:al}, we report the analysis of the variation of the parameter $c_{al}$ with $c_{\rho} = 0.5$.
We do not analyse the case $c_{\rho} = 1.5$ since, as observed in the L-shape test, for this case the variation of $c_{al}$ proves to be irrelevant.
Moreover, we present only the case for VEM order $k=3$, highlighting that the results for $k=1$ and $k=2$ are equivalent.
Figure~\ref{fig:frac86:est:al:dof:05} displays the convergence of the relative a-posteriori error estimator $\eta_{rel}$ of Equation~\eqref{eq:eta:rel} and in Figure~\ref{fig:frac86:est:al:slope:05} we report the convergence rate $\alpha$ at iteration $m>4$, for $c_{al} \in \{0.0, 0.5, 1.0\}$.
As observed in the L-shape test, a value of $c_{al}= 1.0$ facilitates reaching the optimal rate faster than the other values.
Once again, the rates for the other $c_{al}$ values are sub-optimal due to the poor quality of the mesh during the refinement process.
To support this observation, in Figures~\ref{fig:frac86:est:al:q:05}-\ref{fig:frac86:est:al:q:10} we present the statistics for each iteration $m$, showing the evolution of the parameter $\rho_E$ of Equation~\eqref{eq:rhoE}, and the geometric quantities $r_E$ and $h_E$.
We measure the first and third quartile values (q1-q3 areas), and the average values (avg curves).
The measurements confirm that when the length of the minimum edge of the mesh elements ($h_E$) is small compared to the inner radius ($r_E$), i.e. $\rho_E = h_E$, the VEM optimal convergence is not fulfilled.

Setting the optimal value of $c_{al}$ as $1.0$, Figure~\ref{fig:frac86:est:dof_est:05} and Figure~\ref{fig:frac86:est:dof_slope:05} depict the convergence curves for $k \in \{1, 2, 3\}$.
The graphs exhibit the successful convergence of the scheme to the optimal rates $\alpha$ for all orders.
The optimal rate is achieved only after a number of mesh steps, increasing with higher VEM orders: $0.5 \cdot 10^5$, $10^5$, and $0.5 \cdot 10^6$ DOFs for $k=1$, $k=2$, and $k=3$, respectively.
This behaviour is again closely linked to the poor quality of the initial DFN mesh, $\mathcal{T}_{\Omega}^0$, illustrated in Figures~\ref{fig:frac86:m_rhoE:O1:05}-\ref{fig:frac86:m_rhoE:O2:05}.
These figures present for the orders $k = 1$ and $k = 2$ the same statistical analysis of quantities $h_E$, $r_E$, and $\rho_E$ addressed for the case $k=3$ in Figure~\ref{fig:frac86:est:al:q:10}.
When convergence rates $\alpha$ are  sub-optimal (DOFs lower than $10^5$) the average curve of $\rho_E$ does not coincide with the average curve $r_E$.
This discrepancy is more pronounced in real DFN applications where the high randomness of fractures leads to a poor quality discretization. 
Thus, a robust iterative schema, such as the one presented, is crucial to obtain a good numerical solution with the minimum number of DOFs in the discretization.

In conclusion, we assess the computational impact of the proposed Algorithm~\ref{alg:ref} on the global adaptive scheme~\eqref{eq:schema}.
Figure~\ref{fig:frac86:ref:cells:05} illustrates, for each refinement iteration $m$, the comparison between the number of marked cells $\mathcal{M}^m$ to refine ($\mathcal{T}_{\Omega}^{ToRef}$) and the number of refined cells by the algorithm ($\mathcal{T}_{\Omega}^{Ref}$).
Notably, $\mathcal{T}_{\Omega}^{Ref}$ is higher or equal to $\mathcal{T}_{\Omega}^{ToRef}$ due to \emph{extension} proposed in Line~\ref{alg:ref:extend} of \texttt{REFINE} Algorithm~\ref{alg:ref}.
The additional marked elements ($\mathcal{T}_{\Omega}^{Ref} \setminus \mathcal{T}_{\Omega}^{ToRef}$ area) are minimal compared to the size of the initially selected elements.
Moreover, the number of new extended elements approaches to zero after the initial steps.
Additionally, Figure~\ref{fig:frac86:ref:time:05} presents the time ratio employed for the \texttt{SOLVE}, the \texttt{MARK}, and the \texttt{REFINE} steps of the adaptive scheme~\eqref{eq:schema}.
As expected, the \texttt{SOLVE} phase dominates ($\approx 76 \%$).
Remarkably, the time allocated for the new \texttt{REFINE} phase aligns well with that measured in \cite{BERRONE2022103770}, despite the new algorithmic additions.
It is noteworthy that both the \texttt{REFINE} and \texttt{MARK} steps maintain linear complexity with respect to the total number of mesh elements $\# \mathcal{T}_{\Omega}$.

\section{Conclusions}
\label{sec:conclusion}
We present a novel refinement algorithm to enhance the quality of a generic polygonal tessellation combined with the virtual element method.
We show the effectiveness of the new approach in highly intricate domain for flow simulations in fractured media.
It is important to emphasize that the outlined strategy is applicable to a wide range of problems where polygonal meshes are beneficial.

This process is governed by two parameters, which control the shape and the quality of the final discretization elements.
Through numerical examples, we prove the efficacy of the refinement process to improve the representation of the solution, against a starting high number of aligned edges and polygonal cells of varying size. 
This resilience is achieved by incorporating a new control mechanism for the number of the aligned edges.
Moreover, the introduction of marked cell extensions in the refinement algorithm allows us to rapidly obtain a nearly triangular or quadrilateral mesh in few iterations.
These innovative contributions do not compromise the computational time complexity of the method, which remains linear concerning the total number of mesh elements.

In conclusion, our numerical tests reveal optimal convergence rates for high VEM orders, even in challenging scenarios such as the real Discrete Fracture Network benchmark test with non-smooth solutions and initial poor mesh quality.
This highlights the versatility of the proposed refinement algorithm across diverse applications.

\section*{Acknowledgments}
The authors are member of the Gruppo Nazionale Calcolo Scientifico-Istituto Nazionale di Alta Matematica (GNCS-INdAM).
This research has been partially supported by INdAM-GNCS Project CUP:E53C23001670001 and by the MUR PRIN project 20204LN5N5\_003.
The authors kindly acknowledge partial financial support provided by the European Union through project Next Generation EU, PRIN 2022 PNRR project P2022BH5CB\_001 CUP:E53D23017950001.

\bibliographystyle{unsrt}
\bibliography{biblio.bib}

\end{document}